\documentclass[12pt,a4paper]{amsart}
\usepackage{graphicx}
\usepackage{amssymb}
\usepackage{amsthm}
\usepackage{amsmath}
\usepackage{pstcol, pst-node}
\usepackage{mathrsfs}
\usepackage{amscd}
\newtheorem{thm}{Theorem}[section]

\newtheorem{cor}[thm]{Corollary}
\newtheorem{prop}[thm]{Proposition}
\newtheorem{lem}[thm]{Lemma}
\newtheorem{rem}[thm]{Remark}
\newtheorem{claim}[thm]{Claim}

\theoremstyle{definition}
\newtheorem{defn}[thm]{Definition}

\newtheorem{prop-def}[thm]{Proposition-Definition}

\theoremstyle{remark}

\newcommand{\lra}{\longrightarrow}

\newcommand{\be}{\begin{equation}}
\newcommand{\bc}{\begin{cor}}
\newcommand{\bt}{\begin{thm}}
\newcommand{\bl}{\begin{lem}}
\newcommand{\bpr}{\begin{prop}}
\newcommand{\br}{\begin{rem}}
\newcommand{\bd}{\begin{defn}}
\newcommand{\ee}{\end{equation}}
\newcommand{\et}{\end{thm}}
\newcommand{\el}{\end{lem}}
\newcommand{\epr}{\end{prop}}
\newcommand{\er}{\end{rem}}
\newcommand{\ed}{\end{defn}}
\newcommand{\ec}{\end{cor}}
\newcommand{\A}{\Bbb{A}}
\newcommand{\R}{\Bbb{R}}
\newcommand{\C}{\Bbb{C}}
\newcommand{\N}{\Bbb{N}}
\newcommand{\Q}{\Bbb{Q}}
\newcommand{\pll}{\parallel}
\newcommand{\CP}{\Bbb{CP}}

\newcommand{\hs}{\hspace*}

\newcommand{\ep}{\epsilon}
\newcommand{\nnn}{\noindent}

\newcommand{\del}{\partial}

\newcommand{\Z}{\Bbb{Z}}
\begin{document}
\title{Disc counting on toric varieties via tropical curves}
\author{Takeo Nishinou}
\thanks{email : nishinou@kusm.kyoto-u.ac.jp}
\address{Department of Mathematics, Kyoto University,
  Kitashirakawa, Sakyo-ku, Kyoto, Japan} 
\begin{abstract}
In this paper, we define two numbers.
One comes from counting tropical curves with a stop
 and the other is the number of holomorphic discs 
 in toric varieties with Lagrangian
 boundary condition.
Both of these curves should satisfy some matching conditions.
We show that these numbers coincide.
These numbers can be considered as Gromov-Witten type invariants for
 holomorphic discs,
 and they have both similarities and differences to the counting numbers
 of closed holomorphic curves.
We study several aspects of them.
\end{abstract}
\maketitle
\tableofcontents
\nnn
$\bold{Notation.}$ 
We always work over $\Z, \Q, \R$ or $\C$.
We use the same notation as in \cite{NS}.
In particular, $N$ is a free abelian group of rank greater than or
 equal to 2.
$N_{\Q}$ is $N \otimes_{\Z} \Q$ and $N_{\R}$ is $N \otimes_{\Z} \R$.
$M$ is $Hom_{\Z}(N, \Z)$.
$D$ is a closed unit disc in $\C$, which is the domain of  
 stable maps.
A toric variety is always regarded as a complex variety,
 sometimes with a symplectic structure induced from a dual polytope
 (this symplectic structure can be singular
 along lower dimensional toric strata,
 but this does not affect our argument).
If $\Bbb G_{\R}(N)$ is the maximal compact torus of the big torus $\Bbb G(N)$,
 a maximal dimensional
 orbit of $\Bbb G_{\R}(N)$ will be identified with a Lagrangian
 fiber of the moment map.  
If $\Xi\subset N_\Q$ is a subset, $L(\Xi)\subset N_\Q$
 denotes the linear subspace spanned by differences $v-w$ for $v,w\in
 \Xi$, and $C(\Xi)\subset N_\Q\times\Q$ is the closure of the convex
 hull of $\Q_{\ge 0}\cdot (\Xi \times\{1\})$. 
In particular, if $A$ is an
 affine subspace then $L(A)\subset N_\Q$ is the associated linear
 space and $LC(A):=L(C(A))\subset N_\Q\times\Q$ is the linear
 closure of $A\times\{1\}$. 
$\mathbf A$ denotes an affine constraint for tropical curves.
$\widetilde{\mathbf A}$ is the union of an affine constraint and stops,
 with a stop considered as a special constraint.
\section{Introduction}
In \cite{NS}, we considered two numbers.
One is the counting number of genus zero tropical curves in $\R^n$
 and the other is the counting number of 
 rational curves in an $n$-dimensional toric variety.
Both curves should satisfy appropriate matching conditions.
The result of \cite{NS} is that these numbers coincide
 and do not depend on the place of matching conditions.

In this paper, we extend the argument of \cite{NS} to the case
 when the curves have a boundary.
This means we count, on the one hand, holomorphic discs in a 
 toric variety whose boundary is placed on a torus fiber of the
 moment map.
 On the other hand, we count genus zero tropical curves 
 with one `stop' (see Figure 1).
\begin{figure}[h]
\includegraphics{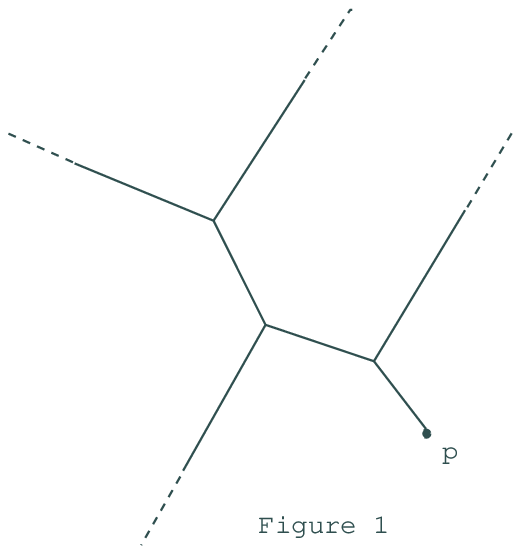}
\end{figure}

These curves should satisfy appropriate matching conditions.
We show that these two numbers coincide.
Moreover, they do not change when the matching condition
 is perturbed slightly.

As a result, we obtain Gromov-Witten type invariants
 for holomorphic discs.
However, these numbers change when the place of the matching condition 
 is changed largely.

In particular, the invariants are well-defined only in the neighbourhood 
 of some degenerations.
We describe these degenerations as toric degenerations (Section 5),
 which should be thought of as the classical limits in the
 context of mirror symmetry.
In this terminology, our invariants for discs may vary as one 
 moves away from the classical limit.
This is a striking difference to the usual Gromov-Witten
 invariants, which are
 invariant under global deformation of the ambient space.

There is one more interesting point to notice.
That is,
 the invariants defined in \cite{M, NS}
 are re-interpretations of Gromov-Witten type invariants
 from tropical geometry.
But the definition of the invariants in this paper relies on tropical 
 geometry and degeneration construction in an essential way.
Some numerical properties of the number of discs are deduced easily via
 tropical geometry,
 while it seems quite difficult to deduce these properties purely
 from complex geometry.\\

The content of this paper is as follows.
In Sections 2 and 3, we give a preparation for tropical curves
 and introduce tropical curves with stops.
In Section 4, we give a preparation for complex curves (with nodes
 and boundary).
In Section 5, we review toric degenerations and holomorphic maps into 
 them.
This is a preparation for Sections 7 and 8.
In Section 6, we treat the family of holomorphic maps from 
 a disc with Lagrangian boundary condition.
The results and techniques are used in Subsection 9.1
 in the proof of the main theorem of this paper.
In Sections 7 and 8, we discuss how one can construct stable maps
 from a disc starting from a tropical curve with a stop.
In Section 9, we give the main theorem, that is,
 the definition of the invariant and give several examples.
In Subsection 9.1, we give the converse of Sections 7 and 8.
Namely, we prove that any stable map from a disc into the toric 
 variety satisfying the required conditions
 is contained in one of the families constructed in Section 8.
This establishes the relation between tropical curves and stable
 maps, and also the invariance of the counting number under perturbation.
In Subsection 9.2, we give examples mainly in two dimensional cases.
We give examples which exhibit the dependence of the invariant on
 the large change of the place of the matching condition,
 the relation to the counting number of closed curves,
 and also give examples where one can calculate the invariant 
 explicitly.
In Subsection 9.3, we give an example in higher dimensional case. 

\nnn {\bf Acknowledgment}
This paper is a continuation of \cite{NS}, and I was benefited
 by a lot of ideas from Bernd Siebert and by discussions with him.
Most importantly, the idea of
 using logarithmic deformation theory in the context
 of tropical geometry is due to him.
He also gave me a lot of useful remarks to this paper.
I would like to express my deep gratitude to him.

\section{Tropical curves with stops}
Here we give the definition of tropical curves and introduce 
 tropical curves with stops, which correspond to holomorphic curves
 with boundary.
Let $\overline \Gamma$ be a weighted, connected finite graph.
Its sets of vertices and edges are denoted $\overline \Gamma^{[0]}$,
 $\overline \Gamma^{[1]}$, and $w_{\overline \Gamma} : 
  \overline \Gamma^{[1]} \to \N \setminus \{ 0 \}$
 is the weight function.
An edge $E \in \overline \Gamma^{[1]}$ has adjacent vertices
 $\del E = \{ V_1, V_2 \}$.
Let $\overline \Gamma^{[0]}_{\infty} \subset \overline \Gamma^{[0]}$
 be the set of one-valent vertices.
We set $\Gamma = \overline \Gamma \setminus \overline \Gamma^{[0]}_{\infty}$.
Non-compact edges of $\Gamma$ are called \emph{unbounded edges}.
Let $\Gamma^{[1]}_{\infty}$ be the set of unbounded edges.
Let $\Gamma^{[0]}, \Gamma^{[1]}, w_{\Gamma}$
 be the sets of vertices and edges of $\Gamma$ and the weight function
 of $\Gamma$, respectively.
\begin{defn}
A \emph{parametrized tropical curve} in $N_{\R}$ is a proper map
 $h : \Gamma \to N_{\R}$ satisfying the following conditions.
\begin{enumerate}
\item[(i)] For every edge $E \subset \Gamma$ the restriction $h \big|_E$
 is an embedding with the image $h(E)$ 
 contained in an affine line with a rational slope.
\item[(ii)] For every vertex $V \in \Gamma^{[0]}$, the following
 \emph{balancing condition} holds.
 Let $E_1, \dots, E_m \in \Gamma^{[1]}$ be the edges adjacent to $V$ and
 let $u_i \in N$ be the primitive integral vector emanating from $h(V)$
 in the direction of $h(E_i)$.
 Then
\begin{equation}
\sum_{j=1}^m w(E_j)u_j = 0.
\end{equation}
\end{enumerate}
\end{defn}
An isomorphism of parametrized tropical curves $h : \Gamma \to N_{\R}$ and 
 $h' : \Gamma' \to N_{\R}$ is a homeomorphism $\Phi : \Gamma \to \Gamma'$
 respecting the weights and $h = h' \circ \Phi$.
A \emph{tropical curve} is an isomorphism class of parametrized 
 tropical curves.
(We often do not distinguish parametrized tropical curves and
 tropical curves.)

The set of \emph{flags} of $\Gamma$ is $F\Gamma = \{(V, E) \big|
     V \in \del E \}$.
By (i) of the definition we have a map
 $u : F\Gamma \to N$ sending a flag $(V, E)$
 to the primitive integral vector $u_{(V, E)} \in N$
 emanating from $V$ in the direction of $h(E)$.

Let $\overline{\Gamma}^{[0]}_{\infty} = \{ q_1, \cdots, q_l, 
 a_1, \cdots, a_m \}$ be the set of
 one-valent vertices.
Let $\Gamma_s = \overline{\Gamma} \setminus \{ q_1, \cdots, q_l \}$.
Let $E_1, \cdots, E_m$ be the edges emanating from $a_1, \cdots, a_m$.
\begin{defn}
A \emph{tropical curve with stops} is the isomorphism class of the proper map
 $h : \Gamma_s \to N_{\R}$ which satisfies
 the conditions of Definition 2.1, with the isomorphism between two 
 maps defined in an obvious way.
We call each $E_i$ a \emph{stopping edge} and $a_i$ a \emph{stop}.
Let $\Gamma^{[1]}_{s, stop}$
 be the set of stopping edges and $\Gamma^{[1]}_{s, \infty}$
 be the set of unbounded
 edges. 
$\Gamma_{s, \infty}^{[1]} = \Gamma_{\infty}^{[1]}
 \setminus \Gamma_{s, stop}^{[1]}$.

We denote by $\Gamma_s^{[0]}$, $\Gamma_s^{[1]}$ the set of vertices
 and edges of $\Gamma_s$.
Let $\Gamma^{[0]}_{stop}$ be the set of stops.
$\Gamma^{[0]}_s =  \Gamma^{[0]} \sqcup \Gamma^{[0]}_{stop}$.
There is a canonical one-to-one correspondence between $\Gamma^{[1]}_s$
 and $\Gamma^{[1]}$.
\end{defn}
An \emph{l-marked} tropical curve (tropical curve with stops)
 is a tropical curve $h : \Gamma \to N_{\R}$
 (resp. $\Gamma_s \to N_{\R}$) together with a choice
 of $l$ distinct edges $\bold{E} = (E_1, \dots, E_l) \in (\Gamma^{[1]})^l$
 (resp. $(\Gamma_{s}^{[1]} \setminus \Gamma_{s, stop}^{[1]})^l$).

The \emph{type} of an $l$-marked tropical curve $(\Gamma, \bold{E}, h)$
 (tropical curve with stops $(\Gamma_s, \bold{E}, h)$)
 is the marked graph $(\Gamma, \bold{E})$ (resp. $(\Gamma_s, \bold{E})$)
 together with the map $u : F\Gamma \to N$ (resp. $u : F\Gamma_s \to N$).

The \emph{degree} of a type
 $(\Gamma, \bold{E}, u)$ ($(\Gamma_s, \bold{E}, u)$)
of a tropical curve (resp. tropical curve with stops) is the function 
 $\Delta: N \setminus \{ 0 \} \to \N$ with finite support defined by
\begin{equation}
\Delta(\Gamma, \bold{E}, u)(v) := \sharp \{(V, E) \in \Gamma^{[1]}_{\infty}
       \big| w(E)u_{(V, E)} = v \}
\end{equation}
\begin{equation}
(resp. \Delta(\Gamma_s, \bold{E}, u)(v) := \sharp \{(V, E) \in 
 \Gamma^{[1]}_{s, \infty}
       \big| w(E)u_{(V, E)} = v \})
\end{equation}
\begin{defn}
The \emph{marked total weight} of a tropical curve with stops is
\begin{equation}
w(\Gamma_s, E) = \prod_{E \in  \Gamma_s^{[1]}
 \setminus \Gamma_{s, \infty}^{[1]}}
                w(E) \prod_{i=1}^l w(E_i).
\end{equation}
\end{defn}

\subsection{Examples}
In this subsection, we give basic examples of tropical curves
 with a stop.
Relations to complex curves (in particular, holomorphic discs
 in toric varieties) are indicated.\\

\nnn
1. {\it Half line}.
The most basic example of a tropical curve with a stop
 is the half line with a rational slope.
Let $e_1, \dots, e_n$ be the standard basis of $\Z^n$.
Then if the slope of the half line is $e_i$ or 
 $- (e_1 + \cdots + e_n)$, then it corresponds to a
 holomorphic disc in $\Bbb P^n$ which intersects at one point
 with a toric divisor
 and whose boundary is mapped to a Lagrangian torus,
 which is a fiber of the moment map.
In fact, half lines are images of such holomorphic discs under
 the moment map (precisely speaking, we identify
 $\R^n$ and the interior of the moment polytope
 via the canonical diffeomorphism 
 between them).
Any holomorphic disc in any toric variety
 which intersects the toric divisor just once
 with the boundary mapped to a Lagrangian fiber of the moment map
 has this property.
This fact follows from, for example, from the explicit description of 
 holomorphic discs with Lagrangian boundary condition given in 
 \cite{CO}, Theorem 5.3.\\

\nnn
2. {\it Discs in $\Bbb P^2$.}
Using the explicit presentation in \cite{CO} above, we can draw the moment map 
 image (i.e, the amoeba) of discs in toric varieties.
For example, the amoeba of a disc with Maslov index 4 is as in Figure 2.
The corresponding tropical curve with a stop is also drawn in the right
 figure.

\begin{figure}[h]
\includegraphics{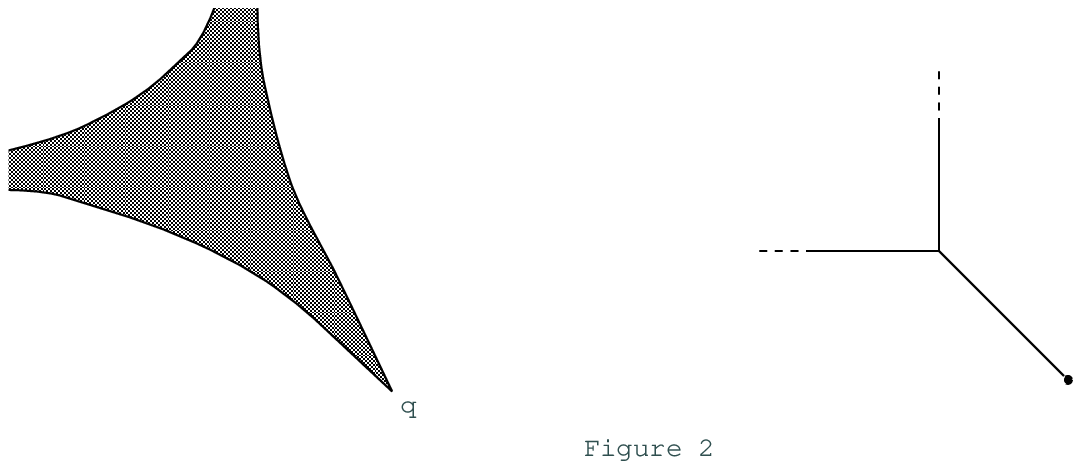}
\end{figure}

One observes that there is a clear resemblance to the amoeba
 or the tropical curve of a line
 in $\Bbb P^2$.
\section{Constraints for tropical curves}
Let $\Gamma_s$ be a weighted graph with stops.
Let $a_1, \dots, a_m$ be the stops of $\Gamma_s$.
Let $\mathbf E = (E_1, \dots, E_l)$ be marking of the edges.
\begin{defn}
For $\bold{d} = (d_1, \dots, d_l) \in \N$ an \emph{affine constraint}
 of codimension $\bold{d}$ is an $l$-tuple
 $\bold{A} = (A_1, \dots, A_l)$ of affine subspaces $A_i \subset N_{\Q}$
 with $dim A_i = n - d_i - 1$.
Set $\widetilde{\bold A} = (p_1, \dots, p_m, A_1, \dots, A_l)$,
 $p_i \in N_{\Q}$.
An $l$-marked tropical curve with stops 
 \emph{matches} the affine constraint 
 $\widetilde{\bold{A}}$ if 
\begin{equation}
\begin{array}{l}
h(E_i) \cap A_i \neq \emptyset, i = 1, \dots l \\
h(a_i) = p_i, i=1, \dots, m
\end{array}
\end{equation}
\end{defn}
The proofs and results in Section 2 of \cite{NS} hold in the 
 case of tropical curves with \emph{one} stop,
 because the slope of the stopping edge
 is uniquely determined by the rest part of the curve 
 by the balancing condition.
In particular, if we extend the stopping edge to infinity,
 we obtain a usual tropical curve.
Then the stop may be considered as a point constraint for this
 tropical curve.
This gives an $(n-1)$-dimensional condition.
Let $\mathfrak T_{(\Gamma_s, \mathbf E,  u)}$ be the set of tropical curves
 with one stop of type $(\Gamma_s, \mathbf E, u)$.
\begin{prop}(Proposition 2.1 of \cite{NS})
For any $\Delta\in Map(N\setminus\{0\},\N)$ and any $g\in\N$ there
 are only finitely many types of tropical curves with one stop
 of degree $\Delta$
 and genus $g$, that is, the set
\[
 \big\{(\Gamma_s, u) \big|
 \Delta(\Gamma_s, u)=\Delta,\ g(\Gamma_s, u)=g,\,
 \mathfrak T_{(\Gamma_s, \mathbf E, u)}\neq\emptyset\big\}
\]
is finite. 
Here $\Gamma_s$ is a weighted graph with one stop, and
 $u: F\Gamma_s \to N$.
\qed
\end{prop}
\begin{rem}
If there are two or more stops, this proposition need not be true.
For example, one constructs tropical curves with two stops and one 
 unbounded edge with the slopes of the edges
 $(k, 1), (-k, 1), (0, 2)$, $k \in \Bbb Z_{>0}$.
See Figure 3.
So there are obviously infinite number of types.
\begin{figure}
\includegraphics{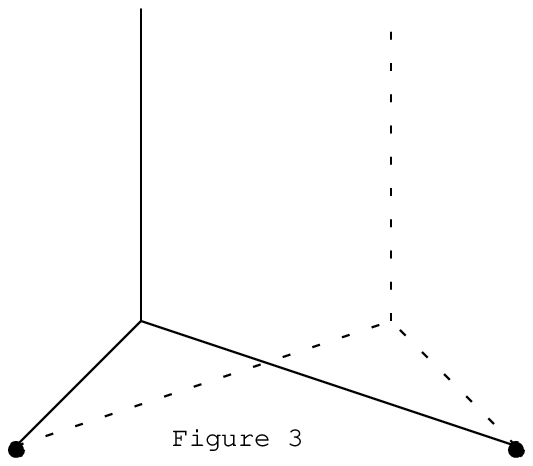}
\end{figure}
\end{rem}
Now we turn to the genus zero case.
As we have mentioned, the degree $\Delta$
 of a tropical curve with a stop uniquely determines the direction 
 $v$ (the primitive integral vector from the vertex which is not the stop)
 of the stopping edge.
Let $\Delta' = \Delta \cup \{ v \} \in Map (\Z^2 \setminus \{ 0 \}, \N)$,
 here $\{ v \} \in Map (\Z^2 \setminus \{ 0 \}, \N)$ is the map
 which takes the value $1$ on $v$ and $0$ otherwise.
\begin{defn}
Let $\Delta\in Map(N\setminus\{0\},\N)$ be a degree and
$e:=|\Delta| = \sharp\Gamma^{[1]}_{s, \infty}$.
An affine constraint $\mathbf{A}=
(A_1,\ldots,A_l)$ of codimension $\mathbf{d}= (d_1,\ldots,d_l)$ 
 and a point $p = A_0 \in N_{\Q}$ is
 \emph{general} for $\Delta$ if
 $\sum_i d_i = (e+1) + n - 3 - (n-1) = e - 1$ and if 
 $\widetilde A$ is general for $l+1$-marked rational tropical 
 curves of degree $\Delta'$.
That is, if
 $(\Gamma_s,\mathbf{E'},h')$ is an $l+1$-marked tropical curve genus $0$ 
 of degree $\Delta'$
 and matching $\widetilde{\mathbf{A}}$, the following holds:
\begin{enumerate}
\item[(i)] $\Gamma$ is trivalent at $\Gamma^{[0]}$.
\item[(ii)] $h'(\Gamma^{[0]})\cap \bigcup_{i=0}^l A_i=\emptyset$.
\item[(iii)] $h'$ is injective for $n>2$. For $n=2$ it is at least
injective on the subset of vertices, and all fibers are finite.
\end{enumerate}
Otherwise it is called \emph{non-general}.\\

The proof of the following result is the same
 for the case with a stop.
\end{defn}
\begin{prop}(Proposition 2.4 of \cite{NS})
Let $\Delta\in Map(N\setminus \{0\},\N)$ and  let $\mathbf A$ be an
 affine constraint of codimension $\mathbf d=(d_1,\ldots,d_l)\in
 \N^l$ with $\sum_i d_i=e-1$. 
Let $\widetilde{\bold A} = \{A_0 = p, A_1, \dots, A_l \}$.
Denote by
 $\widetilde{\A}:=\textstyle{\prod_{i=0}^l} N_\Q/ L(A_i)$ the space of affine
 constraints that are parallel to $\widetilde{\mathbf A}$. Then the subset
\[
 \mathfrak{Z}
  :=\Big\{ \widetilde{\mathbf{A}}'\in
       \widetilde{\A}\,\Big|\, \widetilde{\mathbf{A}}'
 \text{ is non-general for } \Delta\,\Big\}
\]
of $\widetilde{\A}$ is nowhere dense.

Moreover, for any $\widetilde{\mathbf{A}}'\in \widetilde{\A}
 \setminus \mathfrak{Z}$ and any
 $l$-marked type $(\Gamma_s,\mathbf{E},u)$ of genus $0$ and degree
 $\Delta$ there is at most one tropical curve 
 with one stop of
 type $(\Gamma_s,\mathbf{E},u)$ matching $\widetilde{\mathbf{A}}'$. \qed
\end{prop}

\section{Bordered Riemannian surfaces}
We prepare some words for surfaces with boundary.
See also \cite{KL}.
\begin{defn}
Let $A$ and $B$ be nonempty subsets of the upper half plane
 $\Bbb H = \{ z \in \C \big | Im(z) \geq 0 \}$.
A continuous function $f: A \to B$ is \emph{holomorphic} on $A$ if
 it extends to a holomorphic function $\tilde f: U \to \C$,
 where $U$ is an open neighbourhood of $A$ in $\C$.
\end{defn}
\begin{defn}
A \emph{surface} is a Hausdorff, connected, topological space
 $\Sigma$ together with a family
 $\mathcal A = \{(U_i, \phi_i) \big | i \in I \}$
 such that 
 $\{ U_i \big | i \in I \}$
 is an open covering of $\Sigma$ and each map 
 $\phi_i: U_i \to A_i$ is a homeomorphism onto an
 open subset $A_i$ of $\Bbb H$.
The boundary of $\Sigma$ is the set
\begin{equation*}
\del \Sigma = \{ x \in \Sigma \big | \exists i \in I \,\text{such that} \,
 x \in
                    U_i, \phi_i(x) \in \R \} 
\end{equation*} 
The maps
 $\phi_{ij} = \phi_i \circ \phi_j^{-1} :
                   \phi_j (U_i \cap U_j) \to \phi_i (U_i \cap U_j)$
 are surjective homeomorphisms, called
 the transition functions of $\mathcal A$.
$\mathcal A$ is called a \emph{holomorphic atlas}
 if all its transition functions
 are holomorphic. 
\end{defn}
\begin{defn}
A \emph{bordered Riemannian surface} is a compact surface with nonempty
 boundary equipped with the holomorphic structure induced by a 
 holomorphic atlas on it.
A \emph{disc} is a compact surface which is biholomorphic to the unit disc 
 in the complex plane.
\end{defn}
We mainly deal with bordered Riemannian surfaces
 of genus $0$ having one boundary component and some nodes (see below).
\begin{defn}
A \emph{nodal bordered Riemannian surface} is a compact Hausdorff, connected,
 topological space $\Sigma$ with a finite point set $\{ s_1, \dots, s_k \}$
 satisfying the following properties.
Let $\sigma_1, \dots, \sigma_l$ be the connected components of 
 $\Sigma \setminus \{ s_1, \dots, s_k \}$
 ($l$ must be finite due to the following condition at $\{ s_i \}$).
Then $\sigma_i$ is given a structure of a holomorphic surface in the sense of
 Definition 4.2.
At $\{ s_i \}$, there is a neighbourhood
 which is homeomorphic either to
 a neighbourhood of ($0, 0$) $\in \{xy = 0 \}$,
 where ($x, y$) are coordinates on $\C^2$,
 or to a neighbourhood of  
 $(0, 0) \in \{xy = 0 \}/A$, where $A(x, y) = (\overline x, \overline y)$
 is the complex conjugation.
Moreover, the holomorphic structures induced in these neighborhoods
 are compatible with the given holomorphic structures 
 of $\sigma_i$.

The points $\{ s_1, \dots, s_k \}$ are called \emph{nodes}.
\end{defn}
\begin{defn}
A \emph{prestable bordered Riemannian surface}
 is a singular bordered Riemannian 
 surface (see \cite{KaKa} for singularities of 
 Riemannian surfaces) whose singularities are nodes.
A \emph{pointed prestable bordered Riemannian surface} is a prestable 
 bordered Riemannian surface with a collection of points
 $\{ p_i \}_{i = 1}^n$ on it.
No $p_i$ coincides with a node.
A \emph{stable bordered Riemannian surface} is a
 pointed prestable bordered Riemannian surface
 whose automorphism group is finite.
\end{defn}
\begin{defn}
Let $\Sigma$ be a prestable bordered Riemannian surface,
 let ($X, J, \omega$) be a symplectic manifold together with 
 a compatible almost complex structure, and let $L \subset X$
 be a Lagrangian submanifold of $X$.
A \emph{prestable map} $f: (\Sigma, \del \Sigma) \to (X, L)$
 is a continuous map which is $J$-holomorphic 
 in the sense of Definition 4.1, with a straightforward modification.
A \emph{stable map} from a pointed prestable Riemannian surface
 $(\Sigma, \del\Sigma, \{ p_i \})$
 is a prestable map whose automorphism group is finite.
\end{defn}
\br
We can generalize the boundary condition from a Lagrangian submanifold
 to a totally real submanifold.
\er
\subsection{Models for degenerating discs}
In this section, we construct a model for the domain of 
 stable maps from a degenerating family of discs.
Suppose that $C'$ is a tree of rational curves with $2r$ ordinary double
 points $p_1, \cdots, p_r, q_1, \cdots, q_r$.
We assume that $C'$ has
 an anti-holomorphic involution $\sigma$
 whose fixed point locus is a circle and
 which maps
 $p_i$ to $q_i$.
Let $C$ be the quotient of $C'$ by $\sigma$. 
It is known that the universal deformation of $C'$ is parametrized
 by $2r$ dimensional affine plane $\C^{2r}$.
Namely, we can choose coordinates $w_1, \cdots, w_r, z_1, \cdots, z_r$
 on $\C^{2r}$ so that $w_i$ corresponds to the smoothing parameter
 of $p_i$ and $z_i$ corresponds to the smoothing parameter of $q_i$.
Taking a one parameter family $\mathcal X'$ of deformations
 $(f_1(t), \cdots, f_r(t), \overline{f_1(t)}, \cdots, \overline{f_r(t)})$,
 $t \in \C$, here $f_i(t)$ is a holomorphic function with $f(0) = 0$,
 we have a 
 family of rational curves with a family $\Sigma$ 
 of anti-holomorphic involutions
 extending $\sigma$.
The central fiber $X'_0$ is the curve $C'$.
Taking the quotient by $\Sigma$, we have a family of discs 
 degenerating to $C$.
This is our model for the domain of stable maps.

\subsection{Log structures on bordered Riemannian surfaces}
Most of definitions and basic results in log geometry can
 be formulated in analytic setting.
Let $(X, \del X)$ be a complex analytic manifold (or more
 generally a complex space
 in the sense of \cite{KaKa}, Definition 43.2)
 possibly with boundary.
Let $\mathcal O_X^{an}$
 be the sheaf of analytic functions on $X$.
\begin{defn}
A \emph{pre-log structure} on $(X, \del X)$
 is a map of sheaves of monoids
\begin{equation}
\alpha: \mathcal M \to  \mathcal O_X^{an}
\end{equation} 
\end{defn}
\begin{defn}
A pre-log structure $\alpha$ is called a \emph{log structure}
 when the restriction
\begin{equation}
\alpha: \alpha^{-1}((\mathcal O_X^{an})^{\times}) \to
             (\mathcal O_X^{an})^{\times}
\end{equation}
 is an isomorphism.
\end{defn}
A complex analytic manifold (possibly with boundary)
 equipped with a log structure is called a 
 \emph{log complex analytic manifold} and we write it as
 $(X, \del X, \mathcal M)$.

Definitions of a log structure associated to a pre-log structure,
 morphism of log complex analytic manifolds, log smooth morphisms etc.
 generalize to this situation in an obvious manner. 
For more information about log structures, see \cite{KF, KK}.

Recall that, the sheaves of logarithmic derivations
 $\Theta_{(\Sigma, \del \Sigma)}$
 and logarithmic differentials $\Omega^1_{(\Sigma, \del \Sigma)}$ 
 on a prestable bordered Riemannian surface
 $(\Sigma, \del \Sigma)$ are locally free.

\section{Pre-log discs in toric degeneration}
\subsection{Polyhedral decomposition of $N_{\Q}$ and toric degeneration}
First we recall the toric degeneration induced by a
 polyhedral decomposition of $N_{\Q} = \Q^n$, from
 \cite{NS}, Section 3.

A \emph{rational polyhedron} is the solution set in $N_{\Q}$ of
 finitely many linear inequalities $\langle m, \cdot \rangle \geq const.$,
 $m \in M_{\Q}$.
Note that a rational polyhedron may not be bounded.
A \emph{face} of a rational polyhedron is a subset where some of
 the defining inequalities are equalities.
A \emph{vertex} is a zero dimensional face.
When $\Xi$ is a rational polyhedron, we denote by $\mathcal{F}(\Xi)$
 the set of faces of $\Xi$.
A rational polyhedron is called \emph{strongly convex} when 
 it has at least one vertex.
\begin{defn}
A (semi-infinite) \emph{polyhedral decomposition} of $N_{\Q}$
 is a covering $\mathcal P = \{ \Xi \}$
 of $N_{\Q}$ by a finite number of strongly convex rational polyhedra
 satisfying the following properties:\\

(i) If $\Xi \in \mathcal P$ and $\Xi' \subset \Xi$ is a face, then 
 $\Xi' \in \mathcal P$.\\

(ii) If $\Xi, \Xi' \in \mathcal P$, then $\Xi \cap \Xi'$ is a common face of 
 $\Xi$ and $\Xi'$.
\end{defn}
The \emph{asymptotic fan} $\Sigma_{\mathcal P}$ 
 of $\mathcal P$ is defined to be\\
\begin{equation*}
\Sigma_{\mathcal P} := \{ lim_{a \to 0} a\Xi \subset N_{\Q} \big |
 \Xi \in \mathcal P \}.
\end{equation*}
Note that for each $\Xi \in \mathcal P$, the rescaling limit 
 $lim_{a \to 0} a\Xi$ exists in Hausdorff sense and 
 it is either a point (when $\Xi$ is bounded) or a cone (otherwise).
Moreover, we have the following.
\begin{lem}(Lemma 3.2, \cite{NS})
$\Sigma_{\mathcal P}$ is a complete fan. \qed
\end{lem}
Let us define a cone $C(\Xi)$ in $N_{\Q} \times \Q$ by
\begin{equation*}
C(\Xi) = \overline{\{ a \cdot (n, 1) \big | a \geq 0, n \in \Xi \}}.
\end{equation*}
$\overline{(\cdots)}$ means the closure.
This is a strongly convex polyhedral cone.
We define a fan $\widetilde \Sigma_{\mathcal P}$ by
\begin{equation*}
\widetilde \Sigma_{\mathcal P}
 = \{ \sigma \big | \sigma \in \mathcal{F}(C(\Xi)), \Xi \in \mathcal P \}. 
\end{equation*}
We have the following.
\begin{lem}(Lemma 3.3, \cite{NS})
If we identify $N_{\Q}$ with $N_{\Q} \times \{ 0 \}$,
 then
\begin{equation*}
\Sigma_{\mathcal P} = \{ \sigma \cap (N_{\Q} \times \{ 0 \})
                         \big | \sigma \in \widetilde \Sigma_{\mathcal P} \}.
\end{equation*} 
\qed
\end{lem}
The second projection $N_{\Q} \times \Q \to \Q$ defines a map of fans
\begin{equation*}
\widetilde \Sigma_{\mathcal P} \to \{ 0, \Q_{\geq 0} \},
\end{equation*}
 here $\{ 0, \Q_{\geq 0} \}$ is a fan in $\Q$ corresponding to the 
 affine line $\A^1 = \C$.
This induces a toric morphism
\begin{equation*}
\pi: \mathfrak X \to \C.
\end{equation*}
Here $\mathfrak X$ is the toric variety
 associated to the fan $\widetilde{\Sigma_{\mathcal P}}$.
We have the following description of the generic fiber.
\begin{lem}(Lemma 3.4, \cite{NS})
For a closed point $t \in \C \setminus \{ 0 \}$, 
 the fiber $\pi^{-1}(t) \subset \mathfrak X$ with
 the action of $\Bbb G(N) \subset \Bbb G(N \times \Z)$
 is torically isomorphic to $X(\Sigma_{\mathcal P})$. 
Here $X(\Sigma_{\mathcal P})$ is the toric variety associated to the fan
 $\Sigma_{\mathcal P}$.
\end{lem}
We can also describe the central fiber $X_{0}$ using $\mathcal P$.
See \cite{NS}, Section 3.
Here we only note that there is one-to-one correspondence
 between $\mathcal P^{[0]}$ and the components of $X_0$.

As Proposition 3.9 of \cite{NS}, we can fix a polyhedral
 decomposition of $N_{\Q}$ suited to a tropical curve with a stop.
\begin{prop}
Fix a complete fan $\Sigma$.
Let $h : \Gamma_s \to N_{\R}$ be a tropical curve with a stop
 such that the directions of unbounded edges are contained
 in $\Sigma^{[1]}$.
Then there exists a polyhedral decomposition $\mathcal P$ of $N_{\Q}$
 with the asymptotic fan $\Sigma$ such that
 $\cup_{b \in \Gamma^{[\mu]}} h(b) \subset \cup_{\Xi \in \mathcal P^{[\mu]}}
  \Xi, \mu = 0, 1$. 
\end{prop}
\proof
It is easy to see that we can
 graft a tree to the stopping edge so that the resulting graph 
 becomes a tropical curve
 such that the directions of unbounded edges are contained
 in $\Sigma^{[1]}$.
Then apply the construction in the proof of Proposition 3.9 of
 \cite{NS}. \qed\\ 
\subsection{Pre-log discs on the central fiber $X_{0}$}
We extend the definition of pre-log curves to the case
 with boundary.
\begin{defn}
Let $X$ be a toric variety.
A holomorphic curve $C \subset X$ is \emph{torically transverse}
 if it is disjoint from all toric strata of codimension $> 1$.

A stable map $\varphi : C \to X$ is torically transverse
 if $\varphi(C) \subset X$ is torically transverse and
 $\varphi^{-1}(int X) \subset C$ is dense.
\end{defn}
Let $X(\Sigma)$ be a toric variety defined by a fan $\Sigma$
 in $N$.
The toric prime divisors on $X(\Sigma)$ are denoted $D_v$ with $v \in N$
 primitive generators of the rays of $\Sigma$.
\begin{defn}
Let $\varphi : C \to X(\Sigma)$ be a torically transverse holomorphic map.
Here $C$ is a closed surface or a disc.
The \emph{degree} $\delta(\varphi) : N \setminus \{ 0 \} \to \N$
 of $\varphi$ is defined as follows.
For primitive $v \in N$ and $\lambda \in \N$, map $\lambda \cdot v$
 to $0$ if $\Q_{\geq 0}v \notin \Sigma^{[1]}$, 
 and to the number of points of multiplicity
 $\lambda$ in $\varphi^* D_v$ otherwise. 
\end{defn}
\begin{defn}
Let $X_0 = \cup_{v \in \mathcal P^{[0]}} X_{v}$ be the central fiber of the 
 toric degeneration $\mathfrak X \to \C$ defined by an integral
 polyhedral decomposition $\mathcal P$ of $N_{\Q}$.
A \emph{pre-log disc} on $X_0$ is a stable map $\phi: C \to X_0$
 from a prestable bordered Riemannian surface $C$ of
 genus $0$ and $\del C = S^1$
 with the following properties:
\begin{enumerate}
\item[(i)] For any $v$, the projection $C \times_{X_0} X_v \to X_v$
 is a torically transverse stable map.
\item[(ii)] Let $P \in C$ maps to the singular locus of $X_0$.
Then $C$ has a node at $P$, and $\phi$
 maps the two branches $(C', P)$, $(C'', P)$ of $C$ at $P$ to different 
 irreducible components $X_{v'}$, $X_{v''}$.
Moreover, if $w'$ is the intersection index with the toric boundary 
 $D' \subset X_{v'}$ of the restriction $(C', P) \to (X_{v'}, D')$
 and $w''$ accordingly for $(C'', P) \to (X_{v''}, D'')$, 
 then $w' = w''$.
\end{enumerate}
For the purpose of this paper, we require the following additional property.
\begin{enumerate} 
\item[(iii)] Let $C'''$ be the irreducible component with boundary.
Let $X_v$ be the irreducible component of $X_0$ to which $C'''$
 is mapped.
Then the boundary $\del C'''$
 is mapped to a Lagrangian fiber of the moment map
 of $X_v$. 
\end{enumerate}
\end{defn}
In this paper, when we talk about a pre-log disc, it means
 a stable map which satisfies all the properties (i), (ii), and (iii)
 above.
\subsection{Relation between pre-log discs and tropical curves}
A pre-log curve on $X_0$ gives a tropical 
 curve 
 as a dual intersection graph.
See Section 4 of \cite {NS}.
In our case, the dual intersection graph of the pre-log disc
 needs not give a tropical curve with a stop.
But it does give a tropical curve with a stop when the image
 $\phi(C''')$ of the component $C'''$ (see Definition 5.8 (iii))
 intersects the toric boundary of 
 $X_v$ just once.
In particular, this is the case when we treat maximally degenerate
 discs (see Section 7).

\subsection{Boundary condition for the disc}
The boundary condition for the family of degenerating discs is 
 fixed as flows.
Let $p \in \mathcal P^{[0]}$
 be a point which will be the image of the stop of $\Gamma_s$
 and assume $p \in N$.
Let $x \in \mathfrak X$ be a general point.
Let $\Bbb G_{\R}(N)$ be the maximal compact subgroup of
 $\Bbb G(N)$ and let $\Bbb G_{m, p}$ be the one parameter subgroup
 of $\Bbb G(N \times \Z)$ generated by the vector 
 $(p, 1) \in N \times \Z$.
Let $\mathcal L$ be the closure of the orbit
 $(\Bbb G_{\R}(N) \times \Bbb G_{m, p}).x$.
This is the set on which the boundaries of the discs should
 be mapped to.
\subsection{Constraints for the disc}
Constraints for the disc are defined as in \cite{NS}.
Namely, an affine subspace $A\subset N_\Q$ spans the linear
 subspace $LC(A)\subset N_\Q\times\Q$ where the fan
 $\widetilde\Sigma_{\mathcal{P}}$ lives. 
For any general closed point $P$ in the big
 torus of $\mathfrak X$, the closure of the orbit
 $\Bbb{G}(LC(A)\cap (N\times\Z)).P$ defines a subvariety $Z\subset
 \mathfrak X$ projecting onto $\C$. 
Our incidence
 condition is non-trivial intersection with $Z$.
See Section 3 of \cite{NS} for details.

There is one point different from the closed curve case.
Namely, we also put an extra marked point, which corresponds
 to the stop of the tropical curve.
Take $p = A_0 \in N$ and $x \in \mathfrak X$ as in Subsection 5.4.
Then the cone $LC(A_0)$ defines a subgroup
 $\Bbb{G}_{m, p}$.
We set $Z_0$ to be the closure of 
 $\Bbb{G}_{m, p}.x$.
\emph{We will count stable maps from pointed discs
 intersecting the constraints
 $Z_1, \dots, Z_l$ in the interior and $Z_0$ on the boundary}.

\section{Gromov compactness theorem for holomorphic discs}
Here we give a corollary of the result of Ye \cite{Ye}
 in a form suited to our purpose.
Let $\pi : \mathfrak X \to \C$ be a toric degeneration
 with the central fiber $X_0$.
After blowing up along toric subvarieties if necessary,
 we torically embed the family $\pi : \mathfrak X \to \C$
 to $\mathbb P^d \to \C$ for some integer $d$.
This blowing-up process does not affect the subsequent argument.
Let $\{ p_i \}_{i=1}^{\infty}$ be a sequence of points
 in $\C$ converging to the origin. 
Let $\phi_i : D \times \{ p_i \} \to \mathfrak X$
 be a family of prestable maps compatible with the projections
 and assume that the boundary maps to a family $\mathcal L$
 of Lagrangian torus
 fibers as in Subsection 5.4.
This family of Lagrangian torus fibers
 extends to a family of Lagrangian torus fibers in 
 $\Bbb P^d \times \C$ by replacing $G_{\R}(N)$ in Subsection 5.4
 by the maximal compact torus acting on $\Bbb P^d$.

To describe Ye's result, we recall some terminology from his paper
 \cite{Ye}.
Let $\Sigma, \Sigma'$ be prestable bordered Riemannian surfaces.
A continuous surjective map $\varphi : \Sigma \to \Sigma'$
 is called a \emph{node map} if the followings hold.
\begin{enumerate}
\item[(i)] For each node $x \in \Sigma'$, $\varphi^{-1}(x)$ is either a node,
 a simple closed curve in the interior which is disjoint from nodes
 or a simple arc which is disjoint from nodes and has
 its endpoints exactly on the boundary.
\item[(ii)] $\varphi$ is a diffeomorphism away from the curves 
 or the points which are 
 preimages of the nodes.
\end{enumerate}
Let $\{ X, \omega, J \}$ be a symplectic manifold with a compatible almost
 complex structure.
\begin{defn}(The $C^k$ topology on the space of prestable maps)
Let $f : \Sigma \to X$ be a prestable map. 
For each $\ep > 0$, a metric on $\Sigma$ and
 a neighbourhood $\widetilde U$ of the nodes of $\widetilde{\Sigma}$,
 a \emph{neighbourhood} $F$ of $f$ is defined as follows.
A prestable map $\widetilde f : \widetilde{\Sigma} \to X$ belongs to
 $F$ if
\begin{enumerate}
\item[(i)] There is a node map $\varphi : \Sigma \to \widetilde{\Sigma}$.
\item[(ii)] $\pll f - \widetilde f \circ \varphi \pll_{C^k} < \ep$
 on $\Sigma$.
\item[(iii)] $\pll j - \varphi^*\widetilde j \pll_{C^k} < \ep$
 on $\varphi^{-1}(\widetilde{\Sigma} \setminus \widetilde U)$.
\item[(iv)] $|area(f) - area(\widetilde f)| < \ep$.
\end{enumerate}
Here $j, \widetilde j$ denote
 the complex structure of $\Sigma$ and $\widetilde{\Sigma}$.
The norms and areas are defined in terms of the given metrics on
 $X$ and $\Sigma$.
\end{defn}
Let $\mathcal C$ be the space of prestable maps to $\mathfrak X$ 
 (or to $\Bbb P^d \times \C$
 by the above embedding).
The definition above generates a topology in $\mathcal C$.
This is called the $C^k$-topology. 
Then Ye's result gives the following.
\bt\label{properness}
There is a prestable bordered Riemannian surface $C_0$ and 
 a prestable map $\phi_0 : C_0 \to X_0$
 such that there is a subsequence of $\{ \phi_i \}$
 which converges to
 $\phi_0$ in $C^k$-topology. \qed
\et
That is, there are node maps $\varphi_i : D \to C_0$
 such that for any $\widetilde U$ (a neighbourhood of nodes
 in $C_0$), and any $\ep > 0$,
 there exists an integer $M$ such that for all $i \geq M$, 
 the triple $(\phi_i, \phi_0, \varphi_i)$ (these are
 $(f, \widetilde f, \varphi)$ in the above definition)
 satisfies the conditions in Definition 6.1.
Note that $(C_0, \phi_0)$ is not unique. 
\subsection{Regularity of maps}
Let $U$ be a neighbourhood of $0 \in \C$.
Let $U^* = U \setminus \{ 0 \}$.
Let $C^* = D \times U^*$ be the product and $\{ x^*_i \}_{i=1}^{k}$
 be holomorphic sections which do not intersect $\del D$.
Let $x^*_0$ be a real analytic section of $\del D \times U^* \to U^*$.
Let $\phi: (C^*, \{x^*_i \}_{i=0}^k) \to \mathfrak X$
 be a family of stable maps
 over $U^*$ satisfying the constraints given in Subsections 5.4 and 5.5.
Let $\phi_0 : C_0 \to X_0$ be an extension given by Theorem 6.2.
\begin{lem}
$\overline{Im \phi} \cap X_0$ is an analytic subset of dimension 1
 with boundary.
\end{lem}
\proof
As proved in the explicit presentation of \cite{CO} Theorem 5.3, 
 any holomorphic disc with Lagrangian boundary condition in 
 a toric variety can be lifted to $\C^n$ for some $n$ and
 the lift to $\C^m$ are 
 written as $z_i = c_i\cdot \prod^{\mu_j}_{k=1}
 \frac{z-\alpha_{i, k}}{1-\overline{\alpha_{i, k}}z}$
 for $c_i \in \C^*$ and $\alpha_{i, k} \in \{ z \in \C; |z| < 1\}$.
Here $z_i$ is a coordinate for $\C^m$.

Since the incidence varieties $\{ Z_i \}$ and the boundary condition
 $\mathcal L$ are algebraic,
 they give algebraic conditions for the parameters $\{ c_i, \alpha_s \}$
 in the explicit presentation.
The fact that the discs are Fredholm regular (\cite{CO} Theorem 6.1)
 and the genericity of the constraint
 implies the solution $\{ c_i, \alpha_s \}$
 is uniquely determined for the family $\phi$,
 and $c_i, \alpha_s$ are algebraic functions of $t$ and $\overline t$
 ($t$ is the holomorphic coordinate of $U$).
So the solutions
 $\{ \alpha_s \}$ are algebraic functions on $U^*$ with values in 
 $\{ w; |w|<1 \}$.
So if we allow the variable $z$ (the coordinate on the disc)
 to take values in $\Bbb P^1$,
 we obtain a flat family
 $\overline{\phi} : \Bbb P^1 \times U^* \to \mathfrak X \setminus X_0$.
Taking the closure,
 we obtain a flat family of curves in $\mathfrak X$
 over $U$.
In particular, the intersection with the central fiber
 $\overline{Im \overline\phi} \subset X_0$
 is also a curve.
So obviously $\overline{Im \phi} \cap X_0$ is also of dimension 1. \qed\\

We recall a result about transversality of a subset
 in a toric variety.
\begin{prop}(Proposition 6.2 of \cite{NS})
Let $X$ be a toric variety and $W\subset X$ a closed complex
 analytic subset of
 codimension $>c$. We assume that no component of $W$ is
 contained in the toric boundary. 
Then there exists a toric blow-up
 $\Upsilon:\tilde X\to X$ such that the strict transform $\tilde W$ of $W$
 under $\Upsilon$ is disjoint from any toric stratum of dimension $\le c$.
 \qed
\end{prop}
By this proposition and Lemma 6.5,
 we assume that $\phi_0(C_0)$ is torically transverse
 in $X_0$.
\begin{lem}
Consider the family $\phi$ and its extension $\phi_0$ as above.
When the constraint is generic, then no boundary bubble occurs
 at $\phi_0$.
\end{lem}
\proof 
Let $X$ be an $n$-dimensional toric variety
 and $L$ a Lagrangian torus fiber of the moment map.
Let $D$ be a disc and $D \cup D$ be the disjoint union of two discs.
Suppose $\phi_1 : D \to X$ and $\phi_2 : D \cup D \to X$
 be two holomorphic maps with boundaries on $L$,
 and suppose that the images have the same relative homology class
 in $H_2(X, L)$.
Let $\mu$ be the Maslov number of these maps. 
It is equal because the two maps have the images
 with the same homology class.
Then the virtual dimension of the moduli space of maps $\phi_1$
 is $\mu + n - 3$.
On the other hand, the virtual dimension of the moduli space of maps
 $\phi_2$, with the additional condition that the intersection of
 the images of the boundaries of the two discs is not empty, is
 $\mu + 2n - 6 - (n-2) = \mu + n - 4$.
As shown in \cite{CO}, the holomorphic maps
 like $\phi_1$ and $\phi_2$ are regular, which means the dimension of the 
 moduli space is equal to the virtual dimension.
Back to the degeneration, this implies that the boundary bubble
 is a codimension one phenomenon.
It follows that if the constraint is generic,
 the boundary bubble does not occur.\qed
\subsection{Separatedness and properness}
Take $p : (C^*, \{ x^*_i \}) \to U^*$ and
 $\phi : (C^*, \{ x^*_i \}) \to \mathfrak X$
 as above.
By Theorem 6.2, we have at least one extension of $p$ and $\phi$ to the 
 origin.
We want to show that $p$ and $\phi$ extend 
 so that they give holomorphic families.
By Lemma 6.3, we only consider the case with no boundary bubble.
The following is obvious from Lemma 6.3.
\begin{cor}
Any limit bordered Riemannian surface obtained from a sequence
 $\{ p_i \} \subset U$ which appears
 as the image $\phi_0(C_0)$ by Theorem 6.2 does not depend on $\{ p_i \}$
 or the subsequences to be chosen. \qed
\end{cor} 
Moreover, we claim the following.
\begin{lem}
We can assume that
 the limit of the domain of the map, that is, $C_0$ in Theorem 6.2 
 does not depend on the sequence $\{ p_i \}$.
Moreover, the map $\phi_0 : C_0 \to X_0$
 does not depend on $\{ p_i \}$, either.
\end{lem}
\proof
One can choose a basis $\bar t = (t_0, \dots, t_d)$
 of homogeneous coordinates of $\mathbb P^d$ so that 
 they intersect $\overline{Im\phi}$ 
 transversally at nonsingular points on each fiber.
In particular,
 we can assume that
 the hyperplanes $\{ t_i = 0 \}$ do not
 intersect the image of the boundary,
 when we take $U$ small.

Take a totally ramified base change so that the divisors
 $\phi^*(t_i)$ split into sections.
We regard these divisors as additional marked points
 and denote by $\{ s_i^* \}$.
Apply doubling construction of Subsection 4.1 to $C^*$.
Then $C^*$ becomes a real analytic family of Deligne-Mumford stable curves
 $\widetilde C^*$.
So by the properness and separatedness
 of the stack of Deligne-Mumford stable curves,
 it extends uniquely to the origin.
Let $s_i$ be the extension of $s_i^*$.
The anti-holomorphic involution also extends.
We denote the central fiber by $\widetilde C_0$ and the quotient
 by the anti-holomorphic involution by $\overline C_0$.

\begin{claim}
The prestable bordered Riemannian surface
 $C_0$ is obtained by attaching some spheres to $\overline C_0$.
\end{claim}
\proof
Otherwise, there will be marked points $s^*_i \neq s^*_j$
 coming from $\phi^*(t_i)$
 over $U \setminus \{ 0 \}$,
 such that $\overline{\phi(s^*_i)} \cap X_0 \neq
 \overline {\phi(s^*_j)} \cap X_0$, but
 $\overline{s^*_i} \cap C_0 =  \overline{s^*_j} \cap C_0$, a contradiction. 
 \qed
\begin{claim}
On the other hand, if $C_0$ has an unstable irreducible component $S$ and
 $\overline C_0$ does not have a corresponding component, then
 $S$ is mapped to a point by $\phi_0$.
\end{claim}
\proof
If $S'$ is a component of $C_0$ such that $\phi(S')$ is not a point, then
 there will be sections $s^*_i$ as above such that 
 $\overline{s^*_i} \cap int S' \neq \emptyset$.
But such a component will also be a component of $\overline{C_0}$. \qed\\

\nnn
So $\phi_0$ descends to
 $\overline C_0$.
We rewrite $\overline C_0$ as $C_0$ afterwards.   
Let $S$ be a component of $C_0$ which has some marked points from 
 $\{ s_i \}$ on its open part.
Then $\phi_0\big |_S$ is a non-constant map and $\phi_0$ maps
 $s_i$ to $\overline{\phi(s_i^*)} \cap X_0$.
It is clear that we can take $\overline t$ so that these conditions
 determine $\phi_0\big |_S$ uniquely.
If on the other hand $S$ does not have a marked point from $\{ s_i \}$ 
 on its interior, then $S$ is mapped to a point.
So $\phi_0$ is uniquely determined. \qed\\

Summarizing, we have proved the following.
\begin{cor}
Let $\phi : (C^*, \{ x^*_i \})
 \to \mathfrak X \setminus X_0$ be a family of stable maps
 from pointed discs.
Then $(C^*, \{ x^*_i \})$ extends uniquely to 
 a family of prestable pointed discs $(C, \{ x_i \})$ 
 so that $\phi$ extends uniquely to a family of stable maps
 $\Phi : (C, \{ x_i \}) \to \mathfrak X$. \qed
\end{cor}  
\section{Maximally degenerate discs in toric degeneration} 
We extend the notion of maximally degenerate curves
 to the case with boundary 
 (see also Section 5, \cite{NS}).
\begin{defn}\label{def lines}
Let $X$ be a complete toric variety and $Y \subset X$ the toric
boundary. A \emph{line} on $X$ is a non-constant, torically
transverse map $\varphi: \Bbb P^1\to X$ such that $\sharp
\varphi^{-1}(Y)\le 3$.

A \emph{half line} on $X$ is a non-constant, torically transverse map
 $\varphi: D \to X$ such that $\sharp
 \varphi^{-1}(Y) = 1$
 and the boundary $S^1$ is mapped to a Lagrangian torus fiber
 of the moment map. (See Example 2.1.1)
\end{defn}
\begin{defn}\label{maximally degenerate curves}
Let $X_0=\bigcup_{v\in\mathcal P^{[0]}} X_v$ be the central fiber of a toric
 degeneration defined by an integral polyhedral decomposition
 $\mathcal P$. 
A pre-log disc (Definition 5.8) $\varphi: C\to X_0$
 is called \emph{maximally degenerate} if for any
 $v\in\mathcal P^{[0]}$ the
 projection $C\times_{X_0}X_v \to X_v$ is a line or a half line.
For $n=2$, the
 union of two divalent lines intersecting disjoint toric
 divisors is also allowed.
Here a divalent line is a line with $\sharp\varphi^{-1}(Y) = 2$.
\end{defn}
Thus a maximally degenerate curve is a collection of lines
 and a half line, at most
 one ($n=2$: two) for each irreducible component of $X_0$, which match
 in the sense that they glue to a pre-log curve.

\subsection{Counting maximally degenerate curves}
We recall from \cite{NS} the result which gives the number of
 maximally degenerate curves on $X_0$ matching the constraint.
The proof is the same even in the case with a stop.
Let $\mathfrak T_{(\Gamma_s,\mathbf E,u)}(\widetilde{\mathbf{A}})$
 be the set of tropical curves of type
 $(\Gamma_s, \mathbf E, u)$ of fixed degree $\Delta$
 satisfying the constraint
 $\widetilde{\mathbf{A}}$.
\begin{prop}\label{count of maximally degenerate curves}
{\rm 1)}\ \ Let $\Delta\in Map(N\setminus\{0\},\N)$ be a degree and
$\mathbf{\widetilde A}=(A_0,\ldots, A_l)$
 an affine constraint (including the stop) that is general
 for $\Delta$ and assume $l \geq 1$.
If $\mathfrak T_{(\Gamma_s,\mathbf E,u)}(\widetilde{\mathbf{A}})\neq\emptyset$
 for an $l$-marked tree with a stop $(\Gamma_s, \mathbf E)$, then the map
\begin{eqnarray*}\label{lattice map}
Map(\Gamma^{[0]},N)&\lra& \prod_{E\in\Gamma_s^{[1]}
 \setminus \Gamma_\infty^{[1]}} N/\Z u_{(\partial^- E, E)} \times
 \prod_{i=0}^l
    N/ \big(\big(\Q u_{(\partial^- E_i,E_i)}+L(A_i)\big)\cap N\big ),\\
 h &\longmapsto& \big((h(\partial^+E)-h(\partial^-E))_E,
 (h(\partial^-E_i - A_i))_i\big).\nonumber
\end{eqnarray*}
 is an inclusion of lattices of finite index
 $\mathfrak D = \mathfrak D(\Gamma_s, \mathbf E, h, \widetilde{\mathbf A})$.
Here
 $\partial^\pm: \Gamma_s^{[1]}\setminus \Gamma^{[1]}_\infty\to
 \Gamma^{[0]}$ is an arbitrarily chosen orientation of the bounded
 edges, that is, $\partial E= \{\partial^- E,\partial^+ E\}$ for any
 $E\in \Gamma^{[1]}\setminus \Gamma^{[1]}_\infty$. If $E\in
 \Gamma^{[1]}_\infty$ then $\partial^-E$ denotes the unique vertex
 adjacent to $E$ which is not the stop.
\smallskip

\noindent
{\rm 2)}\ \ Assume in addition that $\mathcal P$ is an integral polyhedral
 decomposition of $N_\Q$ with
\[
 h(\Gamma_s^{[\mu]})\subset
 \bigcup_{\Xi\in\mathcal P^{[\mu]}} \Xi,\quad \mu=0,1,\quad\text{ and }\quad
 h(\Gamma_s)\cap A_j\subset \mathcal P^{[0]},\quad j=0,\ldots,l,
\]
and associated toric degeneration
 $\mathfrak X \to \C$
 with the central fiber $X_0$. 
Let $P_j$, $j=0,\ldots,l$ be closed
 points in the big torus of $\mathfrak X$. 
Then $\mathfrak D$
 equals the number of isomorphism classes of maximally degenerate
 curves in $X_0$ which is associated to a
 tropical curve with a stop $\Gamma_s$
 (see Subsection 5.3 or Section 4 of \cite{NS}) and
 intersecting
\[
Z_i:=\overline{\Bbb G (LC(A_i)).P_i}\subset \mathfrak X,
\]
 the closure of the orbit through $P_i$ for the subgroup
 $\Bbb G(LC(A_i))\subset \Bbb G(N\times \Z)$ acting on
 $\mathfrak X$.
\end{prop}
On the other hand, given a degenerating family of stable maps to
 $\mathfrak X \to \C$,
 the central fiber can be made torically transverse.
Namely, the proof of Proposition 6.3 of \cite{NS}
 is valid in our case, too.
So we have the following.
Let $U$ be a neighbourhood of $0$ in $\C$ and
 $U^* = U \setminus \{ 0 \}$.
Let $C^* = D \times U^*$.
\begin{prop}
Let $\mathfrak X \to U$ be a toric degeneration with the central fiber
 $X_0 \subset \mathfrak X$.
Let $(C^* \to \mathfrak X \setminus X_0, (x_0^*, \dots, x_l^*))$
 be a family of
 torically transverse stable maps from discs with $l+1$ marked points 
 defined over $U^*$.
Then possibly after base change, there exists a toric blow-up
 $\widetilde{\mathfrak X} \to \C$ with
 center in $X_0$ with the following property.
$C^*$ extends to a family of stable maps
 $(C \to \widetilde{\mathfrak X}, (x_0, \dots, x_l))$
 over $U$ such that for every irreducible component $\widetilde X_v \subset
 \widetilde X_0$ the projection $C_0 \times_{\widetilde X_0} \widetilde X_v
 \to \widetilde X_v$ is a torically transverse stable map.
Here $C_0$ is the central fiber of the family $C$.  \qed
\end{prop}
\section{Deformation theory}
In this section, we extend the arguments in Section 7
 of \cite{NS} to the case with boundary.
We will work in analytic category.
\subsection{Sheaves on bordered Riemannian surfaces}
\begin{defn}
A \emph{Riemann-Hilbert bundle} $(E, E_{\R})$ of rank $n$
 on a bordered Riemannian surface
 $(\Sigma, \del \Sigma)$ is a smooth complex vector bundle $E$ 
 on $\Sigma$ which is holomorphic (in the sense of Definition 4.1,
 obviously modified to this situation. See \cite{KL}, Definition 3.3.12)
 and a totally real subbundle $E_{\R}$ of $E\big |_{\del \Sigma}$.
\end{defn}
By doubling construction , there is a holomorphic vector 
 bundle $E_{\C} \to \Sigma_{\C}$ with an anti-holomorphic
 involution $\widetilde\sigma: E_{\C} \to E_{\C}$ covering
 $\sigma: \Sigma_{\C} \to \Sigma_{\C}$ such that
 $E\big|_{\Sigma} = E$ and the fixed locus of $\widetilde\sigma$
 is $E_{\R} \to \del\Sigma$ (see \cite{KL} Subsection 3.3).

Let $(\mathcal E, \mathcal E_{\R})$ be the sheaf of sections of 
 $(E, E_{\R})$.
Let $U$ be an open subset of $\Sigma$.
Let $\overline U = U \cup \sigma(U)$ be an open subset of $\Sigma_{\C}$.
Then $(\mathcal E, \mathcal E_{\R})(U)
  = \mathcal E_{\C}(\overline U)^{\widetilde\sigma}$.
On the other hand, $-\widetilde\sigma$
 gives another involution covering $\sigma$  
 and 
 $f = \frac{1}{2}(1+\widetilde\sigma)f + \frac{1}{2}(1-\widetilde\sigma)f$,
 $f \in \mathcal E_{\C}(\overline U)$
 gives the decomposition
 $\mathcal E_{\C}(\overline U)
 = \mathcal E_{\C}(\overline U)^{\widetilde\sigma}
   \oplus \mathcal E_{\C}(\overline U)^{-\widetilde\sigma}$.
The spaces on the right hand side are the spaces of sections of 
 Riemann-Hilbert bundles.

The sheaf cohomology group of $(\mathcal E, \mathcal E_{\R})$,
 obtained as the right derived functors of the global section functor 
 as usual, satisfies the equality
 $dim_{\R}H^q(\Sigma, \del\Sigma; E, E_{\R}) 
   = dim_{\C}H^q(\Sigma_{\C}, E_{\C})$.
See \cite{KL}, Section 3 for more details. 

\subsection{Deforming maximally degenerate curves}
As in \cite{NS} Section 7, unadorned letters denote log spaces or morphisms
 of log spaces, while underlined letters denote the underlying 
 spaces or morphisms.
We put a natural log structure on $\mathfrak X$ coming from the
 toric divisor and let $O_0$ be the standard log point
 (see \cite{NS} Section 7).

We extend the results in \cite{NS}, Section 7 to the case with boundary.
Let $\varphi_0 : (\underline C_0, \{ x_i \}) \to \underline X_0$
 be a torically transverse stable map from a prestable pointed disc. 
\begin{prop}\label{log structures on degenerate curves}
Assume that for every bounded edge
 (containing the stopping edge) $E\subset \Gamma_s^{[1]}$ the
 integral length of $h(E)$ is a multiple of its weight $w(E)$.
Then there are exactly $w(\Gamma_s, \mathbf{E})$ pairwise
 non-isomorphic pairs $[\varphi_0: C_0\to X_0,\{ x_i \}]$
 with underlying stable map isomorphic to
 $(\underline{C_0},\{ x_i \},\underline \varphi_0)$, with $\varphi_0$
 strict wherever $X_0\to O_0$ is strict, and such that the
 compositions $C_0\to X_0\to O_0$ are smooth and integral.
Here $w(\Gamma_s, \mathbf{E})$ is the marked total weight (Definition 2.3).
\end{prop}
The proof of this proposition is almost the same as  
 \cite{NS}, Proposition 7.1 and we omit it.
We only remark that if the stopping edge has weight,
 it contributes to the marked total weight.
This is because the disc has two special point, one from the node and
 the other from the marked point on the boundary, so that it has no
 continuous automorphism.
With this remark at hand,
 the argument of the proof of Proposition 7.1 in \cite{NS}
 applies completely similarly.\\

Let $T_{\mathfrak X/\C},  
 T_{C_0/O_0}, etc.$ be logarithmic tangent bundles
 and let  $T^*_{\mathfrak X/\C}, 
 T^*_{C_0/O_0}, etc.$ be
 logarithmic cotangent bundles.
Let $\Theta_{\mathfrak X/\C}, \Theta_{C_0/O_0},
 \Omega_{\mathfrak X/\C}, \Omega_{C_0/O_0}$
 be the sheaves of sections of 
 these bundles. 
Let $N_{\varphi_0} = \varphi_0^*T_{\mathfrak X/\C}/T_{C_0/O_0}$
 be the logarithmic normal bundle and
 $\mathcal N_{\varphi_0}$
 be the corresponding sheaf of sections.

Let $J_{X(\Sigma)}, J_{C_0}, etc.$ be the complex structures on 
 the spaces $X(\Sigma), C_0, etc.$, respectively.
Let $i_{T^n}: T^n \hookrightarrow X_0$ be the Lagrangian
 torus fiber on which the boundary $S^1 \subset C_0$
 is mapped by $\varphi_0$.
We think $T^n$ as a locally ringed space associated to the 
 underlying real analytic manifold.
Let $\mathcal O_{T^n}$ be the structure sheaf.
$\mathcal O_{T^n}$ is a totally real subsheaf
 (i.e, the sheaf of sections of a totally real subvector-bundle)
 of $i_{T^n}^{*}\Bbb C_{X_0}$, here
 $\Bbb C_{X_0}$ is the trivial line bundle on $X_0$.

Let $O_{\infty} = Spec\C[[\ep]]$ be the formal disc.
Let $\widehat{\mathfrak X} = \mathfrak X \times_{\C} O_{\infty}$.
Then $\Theta_{\widehat{\mathfrak X}/O_{\infty}}
 \simeq N \otimes_{\Z} \mathcal O_{\widehat{\mathfrak X}}$
 and
 $i^{*}_{T^n}(T_{\widehat{\mathfrak X}/O_{\infty}})$
 has a totally real subsheaf 
 $N \otimes_{\Z} i \mathcal O_{T^n}$.
We identify this with the sheaf of sections of the tangent bundle
 $T_{T^n}$ of $T^n$.
Similarly, we identify 
 $N \otimes_{\Z} \mathcal O_{T^n}$
 with the sheaf of sections of the normal bundle
 $J_{X_0}
   T_{T^n}$ of $T^n$
 in $X_0$ (or we may take this as the definition of the normal sheaf).
This induces the splitting
 $i_{T^n}^*T_{\widehat{\mathfrak X}/O_{\infty}}
  \simeq T_{T^n}
        \oplus J_{X_0}T_{T^n}$
 as real vector bundles.
The tangent bundle of $C_0$ also splits on the boundary:
 $i_{S^1}^*T_{C_0} = T_{S^1} \oplus J_{C_0}T_{S^1}$,
 here $i_{S^1} : S^1 \hookrightarrow C_0$ is the inclusion.
$(\mathcal N_{\varphi_0},
 \varphi_0^*(\Theta_{T^n})/\Theta_{S^1})$
 and
 $(\mathcal N_{\varphi_0},
 \varphi_0^*
   (J_{X_0}
     \Theta_{T^n})/J_{C_0}\Theta_{S^1})$
 are (sheaves of sections of) Riemann-Hilbert bundles on $(C_0, S^1)$.
Here $\varphi_0^*(\Theta_{T^n})
 = \varphi^{-1}_{0}\big |_{S^1 \times O_0}
    \Theta_{T^n}
 \otimes_{\varphi^{-1}_{0}\big |_{S^1 \times O_0}\mathcal O_{T^n}}
 \mathcal O_{S^1} $
 and we regard $S^1$ as a locally ringed space associated to the
 underlying real analytic manifold.

\begin{lem}\label{deformation theory}
Let $[\varphi_{k-1}: C_{k-1}/O_{k-1}\to \mathfrak X]$ be a lift of
$[\varphi_0:C_0/O_0\to \mathfrak X]$.
Then the set
 of isomorphism classes of lifts $[\varphi_k: C_k/O_k\to \mathfrak X]$
 restricting to $\varphi_{k-1}$ over $O_{k-1}$ is a torsor under
 $H^0(C_0, \del C_0; N_{\varphi_0}, \varphi_0^*(T_{T^n})/T_{S^1})
   \oplus H^0(C_0,\del C_0; \varphi_0^*T_{\mathfrak X/\C},
            \varphi_0^*(J_{X_0}T_{T^n}))$.
\end{lem}
\proof
Let us assume that we have constructed a deformation of 
 $\varphi_0$ up to order $k-1$.
Let $T^n_{(k-1)} \to O_{k-1, \Bbb R}$ be the family of Lagrangian
 torus fibers on which the boundary
 $S^1 \times O_{k-1} \subset C_{k-1}$ is mapped by $\varphi_{k-1}$. 
Here $O_{k, \Bbb R}$ is the locally ringed space structure
 on $O_{k} = Spec \Bbb C[\ep]/(\ep^{k+1})$ whose structure sheaf
 is given by the real sections of that of $O_{k}$.
As in the proof of Lemma 7.2 of \cite{NS}, the extension
 $C_k/O_k$ of $C_{k-1}/O_{k-1}$ exists and is classified by the
 following extension sequence: 
\begin{equation*}
0 \to \mathcal{O}_{C_0} \to \mathcal{E} \to \Omega^1_{C_{k-1}/O_{k-1}} \to 0.
\end{equation*}
The sheaf $\mathcal E$ will be identified with the cotangent sheaf
 $\Omega^1_{C_k/O_k}$ of $C_k$ tensored by $\mathcal O_{C_{k-1}}$.
However, not every extension constructed in this way is relevant to us.
Namely, the curve should be compatible with the construction of
 Subsection 4.1, that is, the curve should be obtained as the quotient 
 of an anti-holomorphic involution which fixes the boundary $S^1$.
So we should only take those extensions obtained from the following 
 extension sequence of Riemann-Hilbert bundles: 
\begin{equation*}
0 \to (\mathcal{O}_{C_0}, \mathcal O_{S^1})
 \to \mathcal{E}
 \to (\Omega^1_{C_{k-1}/O_{k-1}},
 \Omega^1_{S^1 \times O_{k-1, \Bbb R}/O_{k-1, \Bbb R}}) \to 0.
\end{equation*}
%

The existence of an extension of the map
 $\varphi_{k-1}: C_{k-1} \to \mathfrak X$ to $C_k$ is proved 
 as in the proof of Lemma 7.2 of \cite{NS}, 
 replacing the sheaf cohomology group
 $H^i(C_k, \varphi^*_{k-1}\Theta_{X/\Bbb A^1})$ in \cite{NS} by
 $H^i(C_k, \del C_k; \varphi^*_{k-1}\Theta_{\mathfrak X/\C},
 \varphi^*_{k-1}\Theta_{T^n_{(k-1)}/O_{k-1, \Bbb R}})$.
Here
 $\varphi^*_{k-1}\Theta_{T^n_{(k-1)}/O_{k-1, \Bbb R}}
   = \varphi_{k-1}^{-1}\Theta_{T^n_{(k-1)}/O_{k-1, \Bbb R}}
 \otimes_{\varphi_{k-1}^{-1}\mathcal O_{T^n_{(k-1)}}}
 \mathcal O_{S^1 \times O_{k-1, \Bbb R}}$.

These extensions 
 are classified by the 
 factorizations of
 $\varphi^*_{k-1}\Omega^1_{\widehat{\mathfrak X}/O_{\infty}} \to 
   (\Omega^1_{C_{k-1}/O_{k-1}},
 \Omega^1_{S^1 \times O_{k-1, \Bbb R}/O_{k-1, \Bbb R}})$
 over
 $\mathcal E$, that is,
 $\varphi^*_{k-1}\Omega^1_{\widehat{\mathfrak X}/O_{\infty}}
 \to \mathcal E \to (\Omega^1_{C_{k-1}/O_{k-1}},
 \Omega^1_{S^1 \times O_{k-1, \Bbb R}/O_{k-1, \Bbb R}})$. 
Here the map $\varphi^*_{k-1}\Omega^1_{\mathfrak X/O_{\infty}}
    \to (\Omega^1_{C_{k-1}/O_{k-1}},  
 \Omega^1_{S^1 \times O_{k-1, \Bbb R}/O_{k-1, \Bbb R}})$
 is the composition of the standard map 
 $\varphi^*_{k-1}\Omega^1_{\mathfrak X/O_{\infty}}
  \to \Omega^1_{C_{k-1}/O_{k-1}}$ and the projection map
 (see Subsection 8.1.
Namely, it is an isomorphism for an open subset $U \subset C_{k-1}$
 when $U$ does not intersect the boundary.
When $U$ contains a part of the boundary, we can analytically extend
 a section $f \in \Omega^1_{C_{k-1}/O_{k-1}}(U)$ to an open subset 
 containing $U$ in the doubling construction.
Then we can specify the $\widetilde\sigma$-invariant part).

In order to classify the map
 $\varphi^*_{k-1}\Omega^1_{\widehat{\mathfrak X}/O_{\infty}}
 \to \mathcal E$,
 it suffices to classify the maps 
 $(\varphi^{*}_{k-1} \Omega^1_{\widehat{\mathfrak X}/O_{\infty}},
 \varphi^*_{k-1}\Omega^1_{T^n_{(k-1)}/O_{k-1, \Bbb R}})
 \to \mathcal E
 $
 and
 $(\varphi^{*}_{k-1} \Omega^1_{\widehat{\mathfrak X}/O_{\infty}},
 \varphi^*_{k-1}
   J_{\widehat{\mathfrak X}}\Omega^1_{T^n_{(k-1)}/O_{k-1, \Bbb R}})
 \to \mathcal E$. 

The map $\varphi_0: C_0 \to X_0$ satisfies the following 
 boundary condition.
Namely, $\varphi_0$ induces the following commutative diagrams.
\[
\begin{CD}
\varphi_0^*\Omega^1_{\widehat{\mathfrak X}/O_{\infty}}
    @>>> \Omega^1_{C_0} \\
@AAA @AAA\\
(\varphi_0^*\Omega^1_{\widehat{\mathfrak X}/O_{\infty}},
 \varphi^*_0\Omega^1_{T^n}) @>>>
 (\Omega^1_{C_0}, \Omega^1_{S^1})
\end{CD}
\]
\[
\begin{CD}
\varphi_0^*\Omega^1_{\widehat{\mathfrak X}/O_{\infty}}
    @>>> \Omega^1_{C_0} \\
@AAA @AAA\\
(\varphi_0^*\Omega^1_{\widehat{\mathfrak X}/O_{\infty}},
 \varphi^*_0 J_{X_0}\Omega^1_{T^n}) @>>>
 (\Omega^1_{C_0}, J_{C_0}\Omega^1_{S^1})
\end{CD}
\]
The vertical arrows are inclusions.
The second diagram implies that the composition of the maps 
 $(\varphi_0^*\Omega^1_{\widehat{\mathfrak X}/O_{\infty}},
 \varphi^*_0 J_{X_0}\Omega^1_{T^n}) \to \Omega^1_{C_0}
 \to (\Omega^1_{C_0}, \Omega^1_{S^1})$
 is $0$.
The lifts of $\varphi_0$ should also satisfy this property
 (over $O_{k-1}$).
Write $(\varphi^*_{k-1}\Omega^1_{\widehat{\mathfrak X}/O_{\infty}},
 \varphi^*_{k-1}\Omega^1_{T^n_{(k-1)}/O_{k-1, \Bbb R}})
  = \mathcal F$
 and
 $(\varphi^{*}_{k-1}\Omega^1_{\widehat{\mathfrak X}/O_{\infty}},
 \varphi^*_{k-1}J_{\widehat{\mathfrak X}}
 \Omega^1_{T^n_{(k-1)}/O_{k-1, \Bbb R}})
 = \mathcal F_J$.
Then the boundary condition implies that the maps
 $\mathcal F
 \to \Omega^1_{C_{k-1}/O_{k-1}}$ and
 $\mathcal F_J
 \to \Omega^1_{C_{k-1}/O_{k-1}}$  
 split as 
 $\mathcal F
 \to (\Omega^1_{C_{k-1}/O_{k-1}},
 \Omega^1_{S^1 \times O_{k-1, \Bbb R}/O_{k-1, \Bbb R}})
 \to \Omega^1_{C_{k-1}/O_{k-1}}$ and
 $\mathcal F_J
 \to (\Omega^1_{C_{k-1}/O_{k-1}},
 J_{C_{k-1}}\Omega^1_{S^1 \times O_{k-1, \Bbb R}/O_{k-1, \Bbb R}})
  \to \Omega^1_{C_{k-1}/O_{k-1}}$,
 respectively.

Let us consider the first case.
The $k$-th deformation is given by the following commutative diagram.
\[
\begin{CD}
0 @>>> \ker \varphi_{k-1}^* @>>>
         \mathcal F
 @>\varphi^*_{k-1} >>
 (\Omega^1_{C_{k-1}/O_{k-1}},
   \Omega^1_{S^1 \times O_{k-1, \Bbb R}/O_{k-1, \Bbb R}}) @>>>0\\
 @. @V\alpha VV @VVV @VVV\\
 0@>>> (\mathcal O_{C_0}, \mathcal O_{S_1}) @>>> \mathcal E
 @>>>  (\Omega^1_{C_{k-1}/O_{k-1}}, 
         \Omega^1_{S^1 \times O_{k-1, \Bbb R}/O_{k-1, \Bbb R}})
 @>>>0
\end{CD}
\]
The right vertical arrow is the identity.
This diagram is classified up to isomorphism by the map
 $\alpha \in Hom(ker \varphi^*_{k-1}, (\mathcal O_{C_0}, \mathcal O_{S^1})) 
  = Hom(ker \varphi^*_{0}, (\mathcal O_{C_0}, \mathcal O_{S^1})) 
  = H^0(C_0, \del C_0; N_{\varphi_0}, \varphi_0^*(T_{T^n})/T_{S^1})$.

The $k$-th deformation for the second case is given by the similar diagram.
 \[
\begin{CD}
0@>>> \ker \varphi_{k-1}^* @>>>
         \mathcal F_J
@>\varphi^*_{k-1} >>
 0 \in (\Omega^1_{C_{k-1}/O_{k-1}},
 \Omega^1_{S^1 \times O_{k-1, \Bbb R}/O_{k-1, \Bbb R}})
 @>>>0\\
 @. @V\beta VV @VVV @VVV\\
 0@>>> (\mathcal O_{C_0}, \mathcal O_{S^1}) @>>> \mathcal E
 @>>> (\Omega^1_{C_{k-1}/O_{k-1}}, 
 \Omega^1_{S^1 \times O_{k-1, \Bbb R}/O_{k-1, \Bbb R}})
 @>>>0
\end{CD}
\]
This diagram is classified up to isomorphism by the map
 $\beta \in Hom(ker \varphi^*_{k-1}, (\mathcal O_{C_0}, \mathcal O_{S^1}))
 = Hom(ker \varphi^*_{0}, (\mathcal O_{C_0}, \mathcal O_{S^1})) 
 = H^0(C_0, \del C_0;
 \varphi_0^*T_{\mathfrak X/\C}, \varphi_0^*(J_{X_0}T_{T^n}))$.
  \qed\\


Let $\mathbf x_0 = \{ x_0, \dots, x_l \}$
 be the ordered set of marked points.
Let $\mathbf Z = \{ Z_0, \dots, Z_l \}$
 be the ordered set of incident varieties (see Subsection 5.5). 
\begin{prop}\label{deformation with incidence conditions}
Let $[\varphi_{k-1}: C_{k-1}/O_{k-1}\to \mathfrak X,\mathbf x_{k-1}]$
 be a lift
 of $[\varphi_0:C_0/O_0\to \mathfrak X,\mathbf x_0]$ with $\mathbf x_{k-1}$
 factoring over $\mathbf Z$. 
Then up to
 isomorphism there is a unique lift
 $[\varphi_k: C_k/O_k\to \mathfrak X,\mathbf
 x_k]$ to $O_k$ with $\mathbf x_k$ factoring over $\mathbf Z$.
\end{prop}
\proof
The proof goes along the same way as in \cite{NS}, Proposition 7.3.
It suffices to prove the map
\begin{equation}
\begin{array}{l}
H^0(C_0, \del C_0; N_{\varphi_0}, \varphi_0^*(T_{T^n})/T_{S^1})
   \oplus H^0(C_0, \del C_0; \varphi_0^*T_{\mathfrak X/\C},
                \varphi_0^*(J_{X_0}T_{T^n})) \\
 \hs{1in} \to T_{\mathfrak X/\C, \varphi_0(x_0)}/
   (T_{Z_0/\C, \varphi_0(x_0)} + D\varphi_0(T_{S^1, x_0})) \\
 \hs{1.2in}  \oplus
     \prod_{i=1}^{l}T_{\mathfrak X/\C, \varphi_0(x_i)}/
     (T_{Z_i/\C, \varphi_0(x_i)} + D\varphi_0(T_{C_0/O_0, x_i}))
\end{array}
\end{equation}
 is an isomorphism.

The subgroup $H^0(C_0, \del C_0; T_{C_0}, J_{C_0}T_{S^1})$ of
 $H^0(C_0, \del C_0; \varphi_0^*T_{\mathfrak X/\C},
                \varphi_0^*(J_{X_0}T_{T^n}))$ 
 maps to
 $D\varphi_0(J_{C_0}T_{S^1, x_0})$ bijectively,
 and 
\begin{equation*}
\begin{array}{l}
H^0(C_0, \del C_0; \varphi_0^*T_{\mathfrak X/\C},
 \varphi_0^*(J_{X_0}T_{T^n}))
        /H^0(C_0, \del C_0; T_{C_0}, J_{C_0}T_{S^1}) \\
\hs{.1in}      = H^0(C_0, \del C_0; N_{\varphi_0},
  \varphi_0^*(J_{X_0}T_{T^n})/J_{C_0}T_{S^1})
\end{array}
\end{equation*}
So it suffices to prove
 the map
\begin{equation}
\begin{array}{l}
H^0(C_0, \del C_0; N_{\varphi_0}, \varphi_0^*(T_{T^n})/T_{S^1})
   \oplus H^0(C_0, \del C_0; N_{\varphi_0},
 \varphi_0^*(J_{X_0}T_{T^n})/J_{C_0}T_{S^1})\\
 \hs{1in}  \to \prod_{i=0}^{l}T_{\mathfrak X/\C, \varphi_0(x_i)}/
     (T_{Z_i/\C, \varphi_0(x_i)} + D\varphi_0(T_{C_0/O_0, x_i}))
\end{array}
\end{equation}
 is an isomorphism.

We have to interpret
 the left hand side of (9)
 in terms of toric data.

At a vertex $v \in \Gamma_s^{[0]}$, which is not the stop, the map 
 $\varphi^*_0 \Theta_{\mathfrak X/\C} \to \mathcal N_{\varphi_0}$
 induces a canonical surjection
 $N_{\C} = \Gamma(C_v, N \otimes_{\Z} \mathcal O_{C_v})
 \to \Gamma(C_v, \mathcal N_{\varphi_0} \otimes \mathcal O_{C_v})$ 
 as in \cite{NS}.
Here $C_v$ is the rational component corresponding to $v$
 of the maximally degenerate curve. 
At the stop,
 there are natural surjections
 $i_1: N_{\R} = \Gamma(D, \del D; N \otimes_{\Z}(\mathcal O_D, 
 \mathcal O_{S^1})) \to
 \Gamma(D, \del D; \mathcal N_{\varphi_0} \otimes \mathcal O_{D}, 
 \varphi_0^*(J_{X_0}\Theta_{T^n})/J_{C_0}\Theta_{S^1})$
 and $i_2: iN_{\R} = \Gamma(D, \del D; N \otimes_{\Z}(\mathcal O_D, 
 i\mathcal O_{S^1})) \to
 \Gamma(D, \del D; \mathcal N_{\varphi_0} \otimes \mathcal O_{D}, 
 \varphi_0^*(\Theta_{T^n})/\Theta_{S^1})$
Here $\mathcal O_{S^1}$ is the sheaf of real analytic functions on
 $\del D = S^1$.
$(\mathcal O_D, \mathcal O_{S_1})$ is the sheaf of sections of
 a Riemann-Hilbert bundle.

Write $ \Gamma(D, \del D;
  \mathcal N_{\varphi_0} \otimes \mathcal O_{D}, 
   \varphi_0^*(J_{X_0}\Theta_{T^n})/J_{C_0}\Theta_{S^1}) = \Gamma_1$
 and $ \Gamma(D, \del D;
  \mathcal N_{\varphi_0} \otimes \mathcal O_{D}, 
   \varphi_0^*(\Theta_{T^n})/\Theta_{S^1}) = \Gamma_2$.
Let $v$ be the vertex adjacent to the stop. 
Take $h'_v \in N_{\C}$ and a representative 
 $h_1 + h_2 \in N_{\R} \oplus iN_{\R}$ of
 $[h_1] + [h_2] \in \Gamma_1 \oplus \Gamma_2$.
Then $h'_v$ and $[h_1] + [h_2]$
 glue if and only if 
 $h'_v - (h_1 + h_2) \in L(h(E))\otimes \Bbb C$.
Here $E$ us the edge emanating form the stop.
Since $L(h(E))\otimes \Bbb C$ is isomorphic to
 $\ker i_1 \oplus \ker i_2$, the 
 classes $[h_1], [h_2]$ in $\Gamma_1, \Gamma_2$
 are uniquely determined by $h'_v$.

So arguing similarly to \cite{NS}, Proposition 7.3,
 we see that
 the left hand side of (9) is isomorphic to
 $\ker(\Phi: Map(\Gamma^{[0]}, N_{\Bbb C}) \to 
   \prod_{E\in\Gamma_s^{[1]}
 \setminus \Gamma_\infty^{[1]}} N_{\Bbb C}/\C u_{(\partial^- E, E)})$.
Merging $\Phi$ with (9), we obtain a map of Proposition 7.3,
 tensored with $\Bbb C$.
This proves that (9) is an isomorphism.\qed

\section{Counting invariant of discs} 
We will define two numbers by counting two objects.
One is the number of
 tropical curves of genus zero with one stop and given degree $\Delta$
 which match a constraint $\widetilde{\mathbf A}$
 (we count these with weights).
The other is the number of
 holomorphic discs in a toric variety
 which intersect the toric divisor
 in a way specified by the degree $\Delta$ (Definition 5.7), 
 and which match the constraint $\mathbf Z = \{Z_0, \dots, Z_l \}$
 (see Subsection 5.5).
We claim that these numbers coincide.
\subsection{Main theorem}
We will count the number of tropical curves of genus $0$ with one stop
 with the following properties.
\begin{enumerate}
\item[(i)] degree $\Delta$.
\item[(ii)] match the constraint $\mathbf A$
         of codimension $| \Delta | - 1$. 
\item[(iii)] the stop maps to the fixed point 
 $p \in N_{\Q}$.
\end{enumerate}
Tropical curves must be counted with a weight.
It is given by
\begin{equation}
\mathfrak W(\Gamma_s, \mathbf E, h, \widetilde{\mathbf A})
 = w(\Gamma_s, E) \cdot
 \mathfrak D(\Gamma_s, \mathbf E, h, \widetilde{\mathbf A})
  \cdot \prod_{i=1}^{l}\delta_i.
\end{equation}
Here $\delta_i$ is the index of $\Z u_{(\del^-E_i, E_i)} + L(A_i)\cap N$
 in $(\Q u_{(\del^-E_i, E_i)} + L(A_i))\cap N$
 (see remark 5.8 of \cite{NS}).
In general, this counting number depends on the choice of $\mathbf A$ and $p$.
However, we have the following.
\begin{prop}
For any generic constraint $\mathbf A$ and $p$,
 there is a neighbourhood $W$
 of $(\mathbf A, p)$ in $\widetilde{\mathbb A}$ (see Proposition 3.5)
 such that
 the number of tropical curves (counted with the weight
 $\mathfrak W(\Gamma_s, \mathbf E, h, \widetilde{\mathbf A})$) with one stop
 of genus zero, degree $\Delta$ matching the given constraint
 is constant on $W$.
\end{prop}
\proof
Extend the stopping edge of a tropical curve of genus zero,
 degree $\Delta$ with one stop and
 satisfying $\widetilde{\mathbf A} = \{ p, \mathbf A \}$.
We obtain a rational tropical curve satisfying the constraint
 $\{p, \mathbf A \}$. 
By Definition 3.4, $\{p, \mathbf A \}$
 is generic as a constraint for rational tropical curves.
By the proof of \cite{NS}, Proposition 2.4, the type of the tropical curves
 matching $\widetilde{\mathbf A}$
 will not change when we move $\mathbf A$ and $p$ slightly.
This proves the proposition. \qed\\

Let $\mathcal P$ be an integral polyhedral decomposition of $N_{\Q}$
 which contains all the tropical curves with one stop
 which match $\widetilde{\mathbf A}$
 in $\mathcal P^{[1]}$.
(If $\mathcal P$ is not integral, we can make it so by replacing
 $N$ by $\frac{1}{d}N$ for some integer $d$.
This corresponds to a base change in the toric degeneration.)
Let $\Sigma_{\mathcal P}$ be the asymptotic fan.
Let $\pi: \mathfrak X \to \C$ be the toric degeneration
 defined by the fan $\widetilde{\Sigma_{\mathcal P}}$.
\begin{thm}
There is a neighbourhood $U$ of $0 \in \C$
 such that the number of torically transverse
 holomorphic discs of given degree satisfying the 
 constraint $\mathbf Z$ in the fibers of $\pi^{-1}(U \setminus \{0 \})$
 coincides with the weighted number of tropical curves with one stop of 
 genus zero and degree $\Delta$ satisfying the constraint
 $\widetilde{\mathbf A}$.
Here the weight of a tropical curve is given
 by $\mathfrak W(\Gamma_s, \mathbf E, h, \widetilde{\mathbf A})$
 defined above.
\end{thm}
\proof
As in the proof of \cite{NS}, Proposition 8.3, we want to prove that 
 both numbers are equal to the number of maximally degenerate 
 discs satisfying the constraint and equipped with log structures.
We have constructed in the previous section maximally degenerate discs
 from tropical curves with one stop
 and also constructed logarithmic structures on them, and
 the number of them equals the marked total weight
 of the tropical curve.
There are $\mathfrak D$ of maximally degenerate discs
 associated to the tropical curve and
 for each of maximally degenerate curve there are 
 $\prod_{i=1}^l \delta_i$ intersections with the incidence varieties
 $\mathbf Z$. 
We have also proved that each of these logarithmic curves can be deformed
 to give a family of stable maps defined over a neighbourhood
 $U$ of $0 \in \C$.
This gives $\sum \mathfrak W(\Gamma_s, \mathbf E, h, \widetilde{\mathbf A})$
 families of torically transverse stable maps from a disc.
Here the sum is taken over the isomorphism classes of 
 genus 0 $l$-marked tropical curves of degree $\Delta$
 with one stop matching $\widetilde{\mathbf A}$.
So it suffices to prove the following.
\begin{lem}
There is an open neighbourhood $U$ of $0 \in \C $ 
 such that any torically transverse stable map from a disc to $\pi^{-1}(v)$,
 $v \in U \setminus \{ 0 \}$, satisfying the constraint
 is contained in one of the families above.
\end{lem}
\proof
Since the notion of stacks of stable maps from open Riemannian
 surfaces is not established yet enough for our purpose,
 we can not use it, contrary to \cite{NS}.
Let $\phi_i : D \times p_i \to \mathfrak X$,
 $p_i \in \C$, $p_i \to 0$ as $i \to \infty$
 be a sequence of torically transverse
 stable maps each of which is not contained in any of the 
 families above.
Our purpose is to make a family of stable maps
 $\Phi : D \times U^* \to \mathfrak X \setminus X_0$
 from this sequence.
Here $U^* = U \setminus \{ 0 \}$.
Then it extends to a family over $U$ which is torically transverse, 
 by Corollary 6.10 and Proposition 7.4.
Then by the argument of \cite{NS}, Section 8, we deduce a contradiction.

By \cite{CO}, Theorem 6.1, each $\phi_i: D \to \mathfrak X$ over $p_i \in U$
 can be extended to a neighbourhood of $p_i$.
Let $\mathfrak A_i$ be a punctured disc around the origin $0 \in U$
 which contains $p_i$.
We try to extend $\phi_i$ to the family of maps  
 $D \times \mathfrak A_i \to \mathfrak X$ by successively applying
 \cite{CO}, Theorem 6.1.
In general, there will be bubbles which interrupt such extensions.

As in the proof of Lemma 6.3, $\phi_i$s are given by the solutions of
 a system of algebraic equations
 (the unknowns $\alpha_l$ must be in $\{ z \in \C; |z| < 1\}$),
 algebraically parametrized by the coordinate
 $t$ of $U^*$.
The genericity of the constraint implies that this system
 has only a finite number of solutions which gives
 torically transverse disc for generic $t$ and
 this number is uniformly bounded from above in $U^*$.

Moreover, note that the bubbles happen for such $t$ that
\begin{itemize}
\item{the solution has $l \neq m$ such that $\alpha_l = \alpha_m$}
\item{some $\alpha_l$ of the solution has the norm 1}.
\end{itemize}
In particular, for all large $i$, we can assume that there is no bubble
 in the extension process in $\mathfrak A_i$ above.

So we can extend $\phi_i$ like analytic continuation of 
 an analytic function on $\mathfrak A_i$, and thus obtain a covering
 of $\mathfrak A_i$.
This covering is obviously finite due to the finiteness 
 of the solutions mentioned above.

Now after base change if necessary, 
 we have a torically transverse family of stable maps
 $D \times \mathfrak A_i \to \mathfrak X$
 for large $i$.
By corollary 6.10, this extends to the origin.
Then the proof of \cite{NS}, Proposition 8.3 applies and we see that 
 this family is one of the families obtained from the 
 tropical curves explained above.
This contradicts to the assumption that the sequence 
 $D \times \{ p_i \} \to \mathfrak X$ is not contained in these families.
This proves the theorem. \qed\\
\begin{defn}
Let $N_D(\Delta, \mathbf A, p_0)$ be the number determined by Theorem 9.2.
We call this as the \emph{disc counting number}
 of type $(\Delta, \mathbf A, p_0)$.
\end{defn}
On the tropical side, this number is invariant under small perturbation
 of the constraint $\mathbf A$ and $p_0$
 (Proposition 9.1).
Let $O$ be a neighbourhood of $\widetilde{\mathbf A}$ in the space 
 $\widetilde{\Bbb A}$ of 
 constraints.
On the holomorphic disc side, this number is invariant over
 some neighbourhood of the origin
 $V \subset \C$, with the constraint $\mathbf Z'$
 determined by any $\widetilde{\mathbf A'} \in O$.
It is also invariant under small changes of $\mathbf Z$ by the actions of 
 $\Bbb G(N)$, which is the same as changing $P_j$, $j = 0, \dots, l$
 in Proposition 7.3.2.
\subsection{Examples of calculation of the invariant}
\subsubsection{Disc in $\CP^2$ of Maslov index 4.}
This is a very basic example, exhibiting the dependence of the invariant
 on the place of the constraint.
The constraint is given by one point on $\R^2$ (tropical side)
 or on $\CP^2$ (complex curve side).
On the tropical side, it looks like the following picture. 
\begin{figure}[h]
\includegraphics{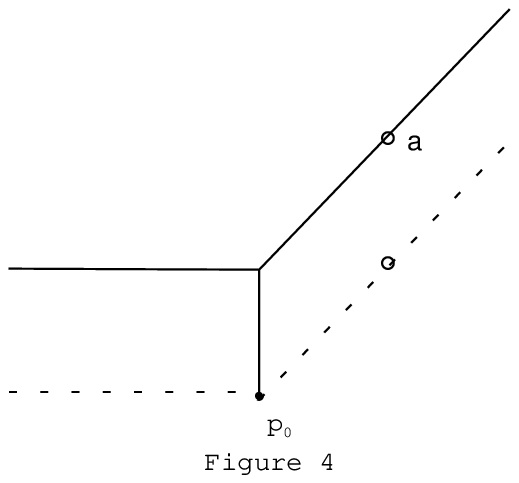}
\end{figure}

The black point $p_0$ is the stop and the circle $a$ is the constraint.
When the circle is sufficiently high above the black point,
 then there is one trivalent tropical curve with one stop.
This means there is a unique smooth disc satisfying the constraint.

When $a$ moves downwards, and intersects the half line of slope $(1, 1)$
 emanating from the stop, then $\widetilde{\mathbf A} = \{ p_0, a \}$
 is not generic.
In this case, the stopping edge degenerates as in the figure.
On the complex curve side, disc bubble occurs on the central fiber
 and the existence of a disc satisfying the constraint
 near the central fiber
 depends on the place of the points $x, p \in \mathfrak X$
 which is used to determine the constraint and the boundary condition
 (Subsections 5.4 and 5.5). 

When $a$ moves downwards further, tropical and holomorphic
 curves satisfying the constraint vanish.
\subsubsection{Count by lattice paths.}
When the constraint is placed in a specific position, 
 we can compute the invariant by Mikhalkin's lattice path count.
Recall that Mikhalkin's lattice path count is applied to count the number of
 tropical curves with the constraint $\{p_1, p_2, \dots \}$
 placed on a line of 
 generic slope, with $d(p_1, p_2) << d(p_2, p_3) << \cdots$
 where $d(p, q)$ is the Euclidean distance between $p, q \in \R^2$.

In our case, we think the stop as one of the constraint,
 and consider the situation when the stop and the constraint are on 
 a line of generic slope.
Let $E$ be the edge of $\Gamma$ containing the stop. 
Let $v \in \Z^2 \simeq N$ be the primitive integral
 vector of $E$ emanating from $p$.
Let $\Delta' = \Delta \cup \{ -v \}$ be the map from
 $\Bbb Z^2 \setminus \{ 0 \}$ to
 $\Bbb N$ which is the union of $\Delta$ and $\{ -v \}
 \in Map(\Bbb Z^2, \Bbb N)$.
Here $\{ -v \}$ means the map
 $\Bbb Z^2 \setminus \{ 0 \} \to \N$ which takes the value 1 on $-v$
 and 0 otherwise.
We count the number of lattice paths on the polytope $P$ dual to the 
 degree $\Delta'$, but with the following condition.
Let $p_1, \cdots, p_r$ be the union of the stop and the constraint,
 and suppose $p = p_i$ is the stop.
The direction $-v$ of the stopping edge determines a dual edge
 $\widehat{-v}$
 of $P$.
Then only the lattice paths whose $i$-th
 path is on $\widehat{-v}$ contributes to
 $N_D(\Delta, \mathbf A, p)$.
So we have the following.
\begin{prop}
In the situation above, $N_D(\Delta, \mathbf A, p_0)$
 is given by the number of lattice paths
 on the polytope dual to the degree $\Delta'$, counted with the weight 
 defined in \cite{M},
 but these paths should satisfy the condition that 
 the i-th path is on $\widehat{-v}$. \qed
\end{prop}

\subsubsection{Relation to closed curves.}
Let $p \in N$ be the image of the stop.
Take $E$ and $\Delta'$ as above.
Extending a tropical curve with a stop at $p$ by extending $E$ 
 to infinity,
 we obtain a tropical curve of genus 0 of degree
 $\Delta'$ without stop.
Let $\widetilde E$ be the extended edge and
 $\widetilde{\Gamma}$ be the extended tropical curve.
The following is clear from the construction.
\begin{prop}
$N_D(\Delta, \mathbf A, p_0) \leq N_{0, \Delta'}^{alg(\mathbf{\widetilde L})}$.
 \qed
\end{prop}
Here $N_{0, \Delta'}^{alg(\mathbf{\widetilde L})}$
 is the number of rational curves
 of degree $\Delta'$ and constraint defined by $\mathbf{\widetilde L}$
 (= linear subspaces of $N_{\R}$ parallel to $\mathbf{\widetilde A}$)
 introduced at Definition 8.2 of \cite{NS}.
Note that this proposition depends heavily on tropical geometry.
It seems difficult to find directly such a bound
 for the number of holomorphic discs.\\

More generally, consider a rational tropical curve $(\Gamma, h)$ in $\R^2$
 and let $E \in \Gamma^{[1]}$
 be the $i$-th marked edge.
Let $r \in E$ be the divalent vertex which is the
 inverse image of the intersection $h(E) \cap p_i$,
 here $p_i$ is the $i$-th constraint.
Splitting $\Gamma$ at $r$, we obtain two tropical curves with a stop
 $(\Gamma_{s, 1}, h|_{\Gamma_{s, 1}})$ and
 $(\Gamma_{s, 2}, h|_{\Gamma_{s, 2}})$.

Conversely, given two tropical curves of genus 0 with one stop whose stopping 
 edges $E_1, E_2$ have the opposite directions, we
 can glue them so that we obtain a rational tropical curve without stop.

The argument in this paper shows we can perform a similar process for 
 holomorphic discs to obtain closed curves,
 though in this case we cannot glue two discs directly.
\begin{figure}[h]
\includegraphics{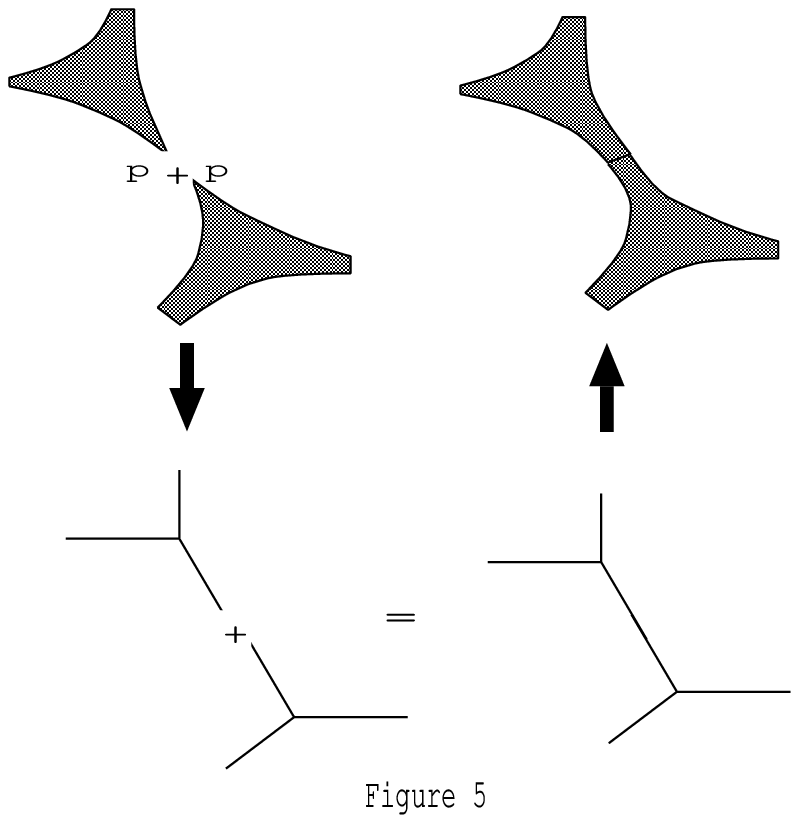}
\end{figure}

In Figure 5, we start from the upper left, which is
 the picture of the amoebas of
 two holomorphic discs $\phi_1, \phi_2$, with boundary on the same
 Lagrangian torus fiber $\pi^{-1}(p)$.
Moreover, their degrees $\delta_1, \delta_2$ sum up to 
 the opposite
 directions:
 $\sum_{v \in N \setminus \{ 0 \}}\delta_1(v)v
   = -\sum_{v \in N \setminus \{ 0 \}}\delta_2(v)v$
 as vectors in $N \setminus \{ 0 \}$.
Apparently we cannot always glue them directly.
But through tropical geometry, we can justify this naive
 gluing process as follows.

Namely, we `tropicalize' them, remembering the logarithmic structures 
 they have (the picture lower left).
In the language of complex curves, this corresponds to 
 considering the maximally degenerate curves (with a logarithmic structure). 
In this tropicalized situation, two discs can be readily glued into 
 one rational tropical curve with a logarithmic structure
 (the picture lower right).

Finally, we smooth back the tropical (or maximally degenerate) curve 
 using the logarithmic structure, and obtain a rational holomorphic 
 curve (the picture upper right).

These rational curves have the  degree which is the sum
 the degrees of the two discs
 and satisfy the constraint which is the sum of the constraints for the
 two discs.
All the rational curves associated to the tropical curve in the
 lower right picture of Figure 5 is obtained in this way
 (over $U$ of Theorem 9.2).

\subsubsection{Constraints in asymptotic position.}
Let $\ell = \{ \ell(a) = p + av | a \in \R_{>0} \}$ be the half line
 emanating from $p$ in the direction $v \in \Q^n$. 
Let $Y_{a, \theta}, a \in \R_{>0}, \theta \in (0, \frac{\pi}{2})$ be the cone 
 defined by
$Y_{a, \theta} = \{ x \in \R^n |
     0 \le \angle(x\ell(a)\ell(a+1)) \le \theta  \}$.  
We consider constraints satisfying the following condition.
\begin{defn}
We say that a constraint $\mathbf A$ satisfies
 the condition $(a, \theta)$ if for each $i$,
 $A_i \cap  Y_{a, \theta} \neq \emptyset$.
\end{defn}
We consider the two dimensional case.
So $\mathbf A$ is a sequence of points.
As above,
 suppose that the stopping edge has direction $v$.
Suppose also that any direction $w$
 of the other unbounded edges (emanating from infinity)
 satisfies the condition
 $\langle v, w \rangle < 0$.
\bpr
When the degree of the tropical curve satisfies this condition,
 there are positive numbers $\eta > 0$ and $\ep << 1$
 such that if $\mathbf A$ satisfies $(a, \theta)$ with
 $a \ge \eta$ and $\theta \le \ep$,
 then the number of genus 0 tropical curves of degree $\Delta$
 which match $\mathbf A$ and satisfy $h(stop) = p$
 is independent of the choice of $\mathbf A$. 
\epr
\proof
Write $\widetilde{\mathbf A} = \mathbf A \cup \{ p \}$.
By the results of \cite{M, NS}, the number of
 rational tropical curves
 of fixed degree which match the constraint $\widetilde{\mathbf A}$
 does not depend on the choice of $\widetilde{\mathbf A}$.
So, to prove the proposition, it suffices to prove the following claim.
\begin{claim}
There is no rational tropical curve of degree
 $\Delta' = \Delta \cup \{ -v \}$ in the notation of 9.2.2
 which matches $\widetilde{\mathbf A}$, with the extra condition that
 an unbounded edge of direction other than $-v$
 or a bounded edge matches $p$.
\end{claim}
Now we prove the claim.
Suppose that there is a tropical curve which matches $\widetilde{\mathbf A}$
 and some unbounded edge $l$ passes $p$,
 whose direction $\vec l$ satisfies $\langle v, \vec l \rangle <0$.
Let $p_1$ be the vertex adjacent to $l$.
Let $H$ be the line through $p$, perpendicular to $v$.
Let $\R^2 = H^+ \cup H^-$ be the decomposition induced by $H$.
By assumption, $\mathbf A$ is contained in $H^+$ or $H^-$. 
Assume $\mathbf A$ is contained in $H^+$.
Then at least one of the edges other than $l$, which emanates from $p_1$
 is contained in $H^-$.
Let $l_1$ be this edge and if $l_1$ is not an unbounded edge,
 let $p_2$ be the vertex adjacent to 
 $l_1$ other than $p_1$.
Again one of the edges emanating from $p_2$ is contained in the half
 plane similar to $H^-$.
Repeating this process, we have to meet an unbounded edge 
 which is contained in $H^-$ (see Figure 6).
By assumption, this edge must have direction $-v$.
Denote this edge by $F$.
\begin{figure}[h]
\includegraphics{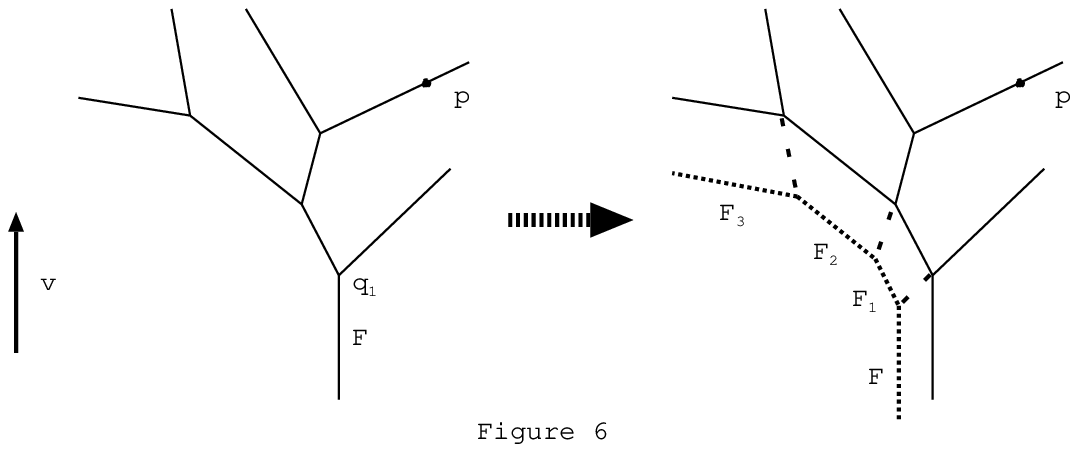}
\end{figure}
\begin{claim}
A tropical curve in $\R^2$ which have a part like this 
 is not rigid. 
That is, there is a continuous deformation keeping the matchings.
\end{claim}
\proof
Suppose that $F$ is placed to the left of $p$.
Let $q_1$ be the vertex adjacent to $F$.
There are two edges, other than $F$, emanating from $q_1$.
One of these edges, say $F_1$, has the direction to the left to $p$.
But such an edge never satisfies $\langle v, w_1 \rangle \le 0$
 where $w_1$ is the direction (from $q$) of the edge $F_1$.
This is because if there is such an edge, it is easy to see that there is 
 an unbounded edge of direction $-v$ other than $F$. 
Similarly, at least one of the edges, say $F_2$ adjacent to $F_1$ has
 the direction to the left to $p$.
If $\theta$ is sufficiently small, $F_2$ and is contained in 
 $\R^2 \setminus Y_{a, \theta}$.
The direction $w_2$ of the edge $F_2$
 will never satisfy $\langle v, w_2 \rangle \leq 0$
 again by the same reason as above.
Continuing this process, we have a family of unmarked edges
 $(F, F_1, \cdots, F_k)$
 with $F$ and $F_k$ unbounded and the others bounded.
We can deform this tropical curve by parallel transforming
 $(F, F_1, \cdots, F_k)$ without changing the incidence conditions
 (see Figure 6).
This means that this tropical curve is not rigid.
When $F$
 is to the right to $p$, do the same by 
  replacing all `left' by `right' in the above argument.
When $F$ is directly below $p$, do the same by using either left or right.
 \qed\\

Since the constraint $\widetilde{\mathbf A}$ is generic,
 there is no non-rigid tropical curve matching $\widetilde{\mathbf A}$,
 by \cite{NS} Proposition 2.4.
So there is no tropical curve of degree $\Delta'$ matching
 $\widetilde{\mathbf A}$ such that an unbounded edge of slope other than 
 $-v$ matches $p$.
For the case when a bounded edge matches $p$, we can do the same
 by using one of the vertices adjacent to this edge.\qed\\
\begin{cor}
Under the above assumption,
 if the degree $\Delta$ and $v$ satisfy the assumption of 
 Proposition 9.12 and $\mathbf A$ whose stop satisfies $(a, \theta)$, then
 $N_D(\Delta, \mathbf A, p) = N_{0, \Delta'}^{alg(\widetilde{\mathbf L})}$.
 \qed
\end{cor} 
When we assume that rational tropical curves of degree $\Delta'$
 satisfying $\widetilde{\mathbf A}$ do not have self-intersection
 (this condition does not depend on the place of $\widetilde{\mathbf A}$),
 we can prove a stronger result by the same argument. 
\begin{prop}
Under the above assumption,
 if $\mathbf A$ and the stop satisfies $(a, \theta)$, then
 $N_D(\Delta, \mathbf A, p) = N_{0, \Delta'}^{alg(\widetilde{\mathbf L})}$.
 \qed
\end{prop}
\begin{rem}
The argument in 9.2.2 and 9.2.3 is valid for any genus tropical curve.
So if we establish a correspondence between genus $g$ tropical curves
 with a stop and genus $g$ holomorphic curves with one boundary component,
 the claims as Proposition 9.5 and corollary 9.11 are also valid.
We will investigate this case,
 together with the cases when there are multiple boundary components
 in a forthcoming paper.
\end{rem}
\subsection{Variation of the invariant}
Although we have considered the cases where we required that the boundary 
 of the disc maps to a particular Lagrangian torus fiber in $X(\Sigma)$,
 we can generalize this.
Namely, instead of $A_0 = p$ in Subsection 5.5, we may take $A_0$ to be an
 $n-d_0-1$ dimensional affine subspace in $N_{\Q}$.
In this case, the constraint should satisfy
 $\sum_{i=0}^l d_i = |\Delta| + n - 2$ and the 
 genericity condition as Definition 3.4.
Then what we count are:
\begin{enumerate}
\item[(tropical side)] genus zero
 tropical curves of given degree with one stop
 matching the constraint $\{A_1, \dots, A_l \}$ and
 the stop is on $A_0$.
\item[(complex side)] stable maps from pointed discs of given degree 
 which meet the subvarieties $\{Z_i \cap X(\Sigma)_t \}$, $i = 1, \dots, l$
 and whose boundary is 
 mapped on a Lagrangian torus contained in $\mathcal L_t$.
Moreover the marked point on the boundary of the disc is mapped to
 $Z_0 = \Bbb G(L(A_0)).P$.
\end{enumerate}
Here $\mathcal L_t$ is the orbit of $\Bbb G_{\R}(N).P$ 
 by the action of $\Bbb G(L(A_0))$, $P \in X(\Sigma)_t$ is a point in the
 big torus. 
$X(\Sigma)_t$ is the fiber over $t \in \C \setminus \{ 0 \}$
 of the family $\pi : \mathfrak X \to \C$. 

A result as Theorem 9.2 holds in this case, too,
 and we can define a disc counting invariant. 
We denote this by $N_D(\Delta, \widetilde{\mathbf A})$
 as before.
Here the weight for the tropical curve should be
 $\mathfrak W(\Gamma_s, \mathbf E, h, \widetilde{\mathbf A})
 = w(\Gamma_s, E) \cdot
 \mathfrak D(\Gamma_s, \mathbf E, h, \widetilde{\mathbf A})
  \cdot \prod_{i=0}^{l}\delta_i$ (there is a contribution 
 $\delta_0$ from $A_0$. 
See Subsection 9.1).
\subsubsection{An example.}
Let us take $p \in \Q^2$, $p_1, \dots, p_l \in \Q^2$,
 a degree function $\Delta_0: \Z^2 \setminus \{ 0 \} \to \N$
 and a vector $v \in \Z^2 \setminus \{ 0 \}$
 so that
 $\{ p, p_i \}$ satisfy the condition $Y_{a, \theta}$
 and the assumption of Proposition 9.8 is also satisfied.
Take a generic decomposition $N_{\Q} = \Q^2 \oplus \Q^{n-2}$
 of $N_{\Q} \simeq \Q^n$
 and let $\pi : N_{\Q} \to \Q^2$ be the projection along this decomposition.
Let us take $A_0 = \pi^{-1}(p)$, $A_i = \pi^{-1}(p_i)$.
Let $\widetilde{\mathbf A} = \{ A_0, \dots, A_l \}$
 and let $\Delta$ be a degree function
 $\Delta : N \setminus \{ 0 \} \to \N$.
Let $w \in N$ be a primitive vector and 
 let $\mathfrak p(w) \in \Q_{\ge 0}.\pi(w)$ be the primitive vector in 
 $\Z^2 \subset \Q^2$, if $\pi(w) \neq 0$.
If there is no $w \in N$ such that $\Delta(w) \neq 0$
 and $\pi(w) = 0$, 
 we can associate to $\Delta$ a degree function
 $\Delta^{\pi}$ on $\Z^2 \subset \Q^2$
 by setting $\Delta^{\pi}(mv') =
 \sum_{\mathfrak p(w) = v'}\Delta(mw)$
 for $v' \in \Z^2$, primitive integral vectors.
Here the sum is taken over primitive $w \in N$.
Take $\Delta$ 
 so that $\Delta$ satisfies $|\Delta| = l + 3 - n$,
 $\Delta^{\pi}$ is defined and satisfies $\Delta^{\pi} = \Delta_0$.
Then by genericity of $\pi$, any genus zero tropical curve
 in $N_{\R}$ with one stop satisfying $\widetilde{\mathbf A}$
 of degree $\Delta$
 will project to a tropical curve on $\Q^2$
 satisfying $\{p, p_i \}$, degree $\Delta_0$.
So by Propositions 9.8 and 9.12, we have the following.
\begin{prop}
Suppose $\Delta_0$, $v$, and $p, p_1, \dots, p_l$ satisfies either
 the assumptions of Proposition 9.8 or 9.12.
Then we have 
 $N_D(\Delta, \widetilde{\mathbf A})
 = N^{alg(\widetilde{\mathbf L})}_{0, \Delta'}$. \qed
\end{prop}

\end{document}